\documentclass[a4paper,11pt]{article}

\usepackage{amsmath,amssymb}

\usepackage{graphicx}
\usepackage{epsfig}
\usepackage{color}

\newtheorem{dfn}{Definition}[section]
\newtheorem{thm}{Theorem}
\newtheorem{lem}{Lemma}[section]
\newtheorem{cor}{Corollary}
\newtheorem{prp}{Proposition}

\newenvironment{proof}{\begin{trivlist}
    \item[\noindent]{\it Proof\:}}{\quad $\square$\end{trivlist}}

\graphicspath{{Figures/}}

\def\R{{\mathbb R}}

\def\Sph{{\mathbb S}}
\def\phi{\varphi}
\def\epsilon{\varepsilon}
\def\l{\ell}
\def\dist{\mathrm{dist}}
\def\D{{\mathcal D}}
\def\P{{\mathcal P}}
\def\Q{{\mathcal Q}}
\def\V{{\mathcal V}}
\def\E{{\mathcal E}}
\def\F{{\mathcal F}}
\def\ext{\mathrm{ext}}
\def\vol{\mathrm{vol}}
\def\per{\mathrm{per}}
\def\area{\mathrm{area}}
\def\lk{\mathrm{lk}\ }
\def\diam{\mathrm{diam}}

\title{Alexandrov's theorem,\\ weighted Delaunay triangulations,\\ and mixed volumes}
\author{Alexander I.~Bobenko\footnote{Institut f\"ur Mathematik, Technische Universit\"at Berlin, Str. des 17. Juni 136, 10623 Berlin, Germany; {\tt bobenko@math.tu-berlin.de}, {\tt izmestiev@math.tu-berlin.de}
\medskip \newline
Research for this article was supported by the DFG Research Unit 565 ``Polyhedral Surfaces''.}  \and Ivan Izmestiev\footnotemark[\value{footnote}]}

\begin{document}

\maketitle

\begin{abstract}
We present a constructive proof of Alexandrov's theorem regarding the existence of a convex polytope with a given metric on the boundary. The polytope is obtained as a result of a certain deformation in the class of generalized convex polytopes with the given boundary. We study the space of generalized convex polytopes and discover a relation with the weighted Delaunay triangulations of polyhedral surfaces. The existence of the deformation follows from the non-degeneracy of the Hessian of the total scalar curvature of a positively curved generalized convex polytope. The latter is shown to be equal to the Hessian of the volume of the dual generalized polyhedron. We prove the non-degeneracy by generalizing the Alexandrov-Fenchel inequality. Our construction of a convex polytope from a given metric is implemented in a computer program.
\end{abstract}

\section{Introduction}
In 1942 A.D.Alexandrov \cite{Al42} proved the following remarkable
theorem:
\begin{thm} \label{thm:Alex}
Let $M$ be a sphere with a convex Euclidean polyhedral metric.
Then there exists a convex polytope $P \subset \R^3$ such that the
boundary of $P$ is isometric to $M$. Besides, $P$ is unique up to
a rigid motion.
\end{thm}
Here is the corresponding definition.

\begin{dfn}
A Euclidean polyhedral metric on a surface $M$ is a metric
structure such that any point $x \in M$ posesses an open
neighborhood $U$ with one of the two following properties. Either
$U$ is isometric to a subset of $\R^2$, or there is an isometry
between $U$ and an open subset of a cone with angle $\alpha \ne
2\pi$ such that $x$ is mapped to the apex. In the first case $x$
is called a regular point, in the second case it is called a
singular point of $M$. The set of singular points is denoted by
$\Sigma$.

If for any $x \in \Sigma$ the angle at $x$ is less than $2\pi$,
then the Euclidean polyhedral metric is said to be convex.
\end{dfn}
For brevity we usually write convex polyhedral metric omitting the
word Euclidean.

The polytope $P$ in Theorem \ref{thm:Alex} may be a plane convex
polygon. In this case the boundary of $P$ consists of two copies
of this polygon identified along the boundary.

The metric on the boundary of a convex polytope is an example of a
convex polyhedral metric on the sphere. The only singular points
are the vertices of the polytope, because any point in the
interior of an edge posesses a Euclidean neighborhood. The
intrinsic metric of the boundary of a polytope does not
distinguish the edges.

In practice, one specifies a polyhedral metric on the sphere by
taking a collection of polygons, identifying them along some pairs
of edges. Alexandrov calls this representation a development, in
analogy with developments of polytopes. By Theorem \ref{thm:Alex},
from any development with cone angles less than $2\pi$ a convex
polytope can be built. However, the edges of this polytope can
differ from the edges of the development. To see this, one can
take a ready-made polytope and cut it into plane pieces by
geodesics different from the edges of the polytope. For a given
development it is a hard problem to determine the edges of the
corresponding polytope.

If a development and the true edges of the polytope are known,
then the problem of constructing the polytope is still
non-trivial, and is studied in \cite{FP} using polynomial
invariants.

The uniquness part of Theorem \ref{thm:Alex} can be proved by a
modification of Cauchy's proof of global rigidity of convex
polytopes. The original existence proof by Alexandrov \cite{Al42},
\cite{Al05} is more involved and non-constructive. Alexandrov in
\cite[pp. 320--321]{Al05} discusses a possibility for a
constructive proof.

Alexandrov's theorem is closely related to the following problem
posed by Hermann Weyl in \cite{Weyl}: \emph{Prove that any
Riemannian metric of positive Gaussian curvature on the sphere can
be realized as the metric of the boundary of a unique convex body
in $\R^3$.} Alexandrov in \cite{Al42} applied Theorem
\ref{thm:Alex} to solve Weyl's problem for the metric of class
$C^2$ through approximation by convex polyhedral metrics. The
usual approach to Weyl's problem through a PDE was realized by
Nirenberg \cite{Nir}.

An interesting remark regarding Weyl's problem was made by
Blaschke and Herglotz in \cite{BH}. Consider all possible
extensions of the given Riemannian metric on the sphere to a
Riemannian metric inside the ball. On the space of such extensions
take the total scalar curvature functional, also known as the
Einstein-Hilbert action. It turns out that its critical points
correspond to Euclidean metrics in the interior of the ball. It is
not clear how to use this observation to solve Weyl's problem
since the functional is not convex.

Since Alexandrov's theorem is a discrete version of Weyl's
problem, one can try to discretize the Blaschke-Herglotz approach.
The discrete analog of the Einstein-Hilbert action is known as the
Regge action. If simplices of a triangulated closed 3-manifold are
equipped with Euclidean structures, then the Regge action
\cite{Reg61} is defined as
\begin{equation} \label{eqn:Reg}
\sum_e \ell_e \kappa_e.
\end{equation}
Here the sum ranges over all edges of the triangulation, $\ell_e$
is the length of the edge $e$, and $\kappa_e$ is the angle
deficit, also called curvature.

This way to discretize the total scalar curvature was also known
to Volkov, a student of Alexandrov. In \cite{Vol68} (translated
and reprinted as \cite[Section 12.1]{Al05}), he applies it to
study the dependence of the extrinsic metric of a convex polytope
on the intrinsic metric on its boundary. Volkov in his PhD of 1955
also gives a new proof of Alexandrov's theorem. It proceeds by
extending the metric to a polyhedral metric in the ball in a
special way and minimizing the sum of edge lengths over all such
extensions. The discrete total scalar curvature functional is not
used in the proof. Volkov's proof can be found in \cite{VP}. See
also Zalgaller's remark in \cite[pp. 504--505]{Al05}.

In this paper we investigate the total scalar curvature functional
(\ref{eqn:Reg}) and give a new proof of Alexandrov's theorem based
on the properties of this functional.

Let us sketch our proof.

A {\em generalized polytope} is a simplicial complex glued from
pyramids over triangles of some geodesic triangulation $T$ of $M$.
Vertices $i \in \V(T)$ of $T$ are singularities of $M$; edges $ij
\in \E(T)$ are geodesic arcs. A generalized polytope can be
described by a couple $(T,r)$, where $r = (r_i)_{i\in \Sigma}$ are
lengths of the side edges of the pyramids. A {\em generalized
convex polytope} is a generalized polytope with all of the
dihedral angles $\theta_{ij}$ at the boundary edges less or equal
to $\pi$. In Section \ref{sec:CGPandDel} we show that for any $r$
there is at most one generalized convex polytope and describe the
space of all generalized convex polytopes with boundary $M$.

In Section \ref{sec:MixVol} we study the total scalar curvature of
generalized convex polytopes:
$$
H(P) = \sum_{i \in \V(T)} r_i \kappa_i + \sum_{ij \in \E(T)}
\ell_{ij} (\pi - \theta_{ij}),
$$
where $\kappa_i$ is the curvature at the $i$-th radial edge, and
$\ell_{ij}$ is the length of the boundary edge joining the
vertices $i$ and $j$. The main result of Section \ref{sec:MixVol}
is the non-degeneracy of the Hessian of $H$ if
\begin{equation} \label{eqn:Cond}
0 < \kappa_i < \delta_i \quad \mbox{ for every } i.
\end{equation}
Here $\delta_i$ is the angle deficit at the $i$-th singularity of
the metric of $M$.

Section \ref{sec:Proof} contains the proof of Alexandrov's
theorem. The idea is to start with a certain generalized convex
polytope over $M$ and deform it into a convex polytope. As the
starting point we take the generalized polytope $(T_D,r)$, where
$T_D$ is the Delaunay triangulation of $M$, and $r_i = R$ for
every $i$, with a sufficiently large $R$. We show that this
generalized polytope is convex and satisfies condition
(\ref{eqn:Cond}). In order to show that any small deformation of
the curvatures $\kappa_i$ can be achieved by a deformation of
radii $r_i$ it suffices to prove the non-degeneracy of the
Jacobian $\left( \frac{\partial \kappa_i}{\partial r_j} \right)$.
Here the functional $H$ comes into play since
$$
\frac{\partial \kappa_i}{\partial r_j} = \frac{\partial^2
H}{\partial r_i \partial r_j}.
$$
We deform $\kappa_i$ by the rule
\begin{equation} \label{eqn:Rule}
\kappa_i(t) = t \cdot \kappa_i,
\end{equation}
where $t$ goes from $1$ to $0$. In particular, condition
(\ref{eqn:Cond}) remains valid during the deformation. We show
that when $r$ is changed according to (\ref{eqn:Rule}), the
generalized convex polytope does not degenerate. In the limit at
$t \to 0$ we get a convex polytope with boundary $M$.

A new proof of Alexandrov's theorem is not the only result of this
paper. Coming across {\em weighted Delaunay triangulations} and
{\em mixed volumes} while studying the total scalar curvature
functional was a big surprise for us. In Section
\ref{sec:CGPandDel} we give a new geometric interpretation of
weighted Delaunay triangulations through generalized convex
polytopes. Also, for a given polyhedral surface with marked
points, we explicitely describe the space of weights of weighted
Delaunay triangulations. Further, we show that a weighted Delaunay
triangulation with given weights can be obtained via a flip
algorithm.

In Section \ref{sec:MixVol} for any generalized convex polytope
$P$ we define its dual generalized convex polyhedron $P^*$. We
find a surprising relation between the total scalar curvature of
$P$ and the volume of $P^*$:
\begin{equation} \label{eqn:PP*}
\frac{\partial^2 H}{\partial r_i \partial r_j}(P) =
\frac{\partial^2 \vol}{\partial h_i \partial h_j}(P^*).
\end{equation}
Here the $h_i$'s are the altitudes of $P^*$, see Subsection
\ref{subsec:GenPolyh}. We prove that the Hessian of $\vol(P^*)$ is
non-degenerate if the curvatures satisfy (\ref{eqn:Cond}). The
proof uses the theory of mixed volumes: we generalize the
classical approach \cite{Schn} to the Alexandrov-Fenchel
inequalities. Note that mixed volumes here play a role similar to
that in the variational proof of Minkowski's theorem, see
\cite[Section 7.2]{Al05}.

\begin{figure}[ht] \label{fig:Examples}
\includegraphics[width=0.26\textwidth,angle=-90]{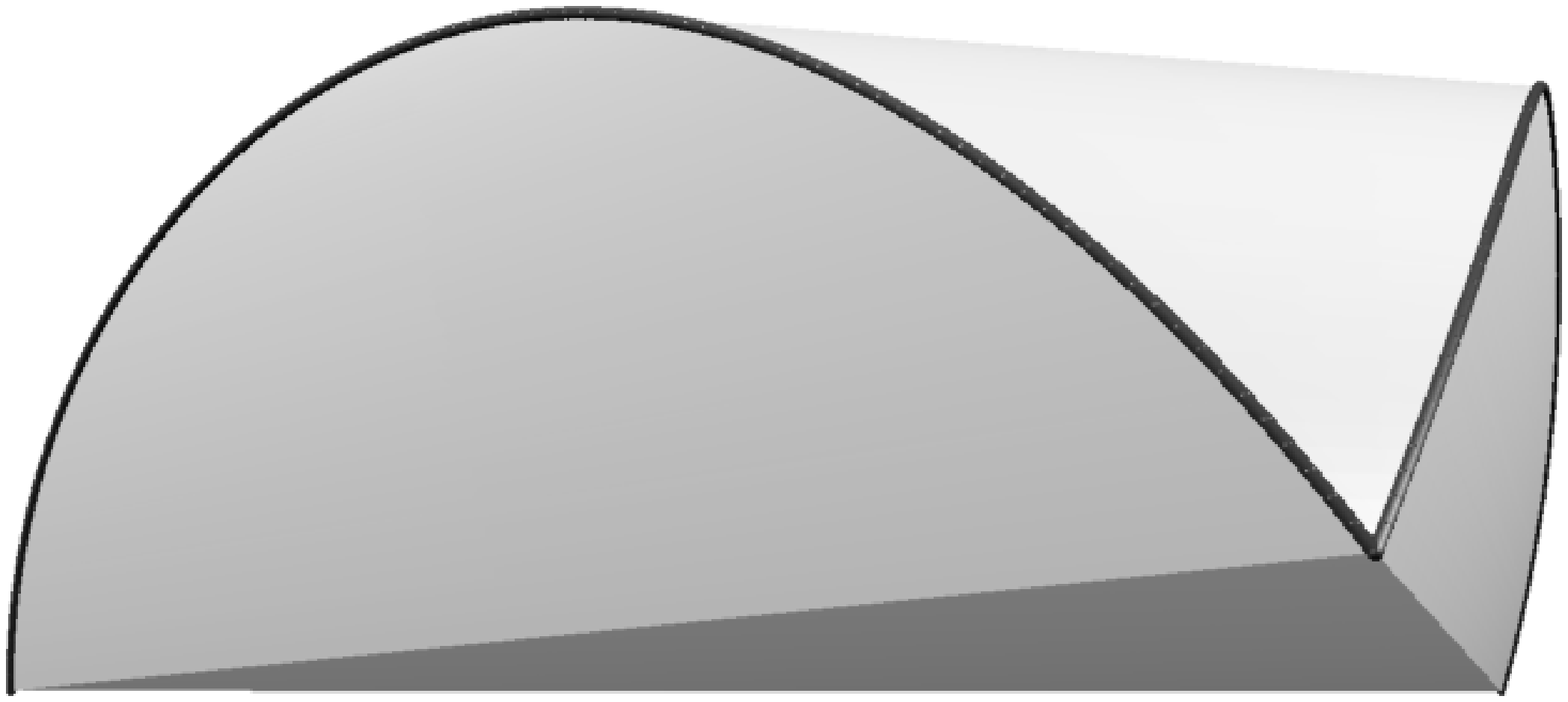} \hfill
\includegraphics[width=0.26\textwidth,angle=-90]{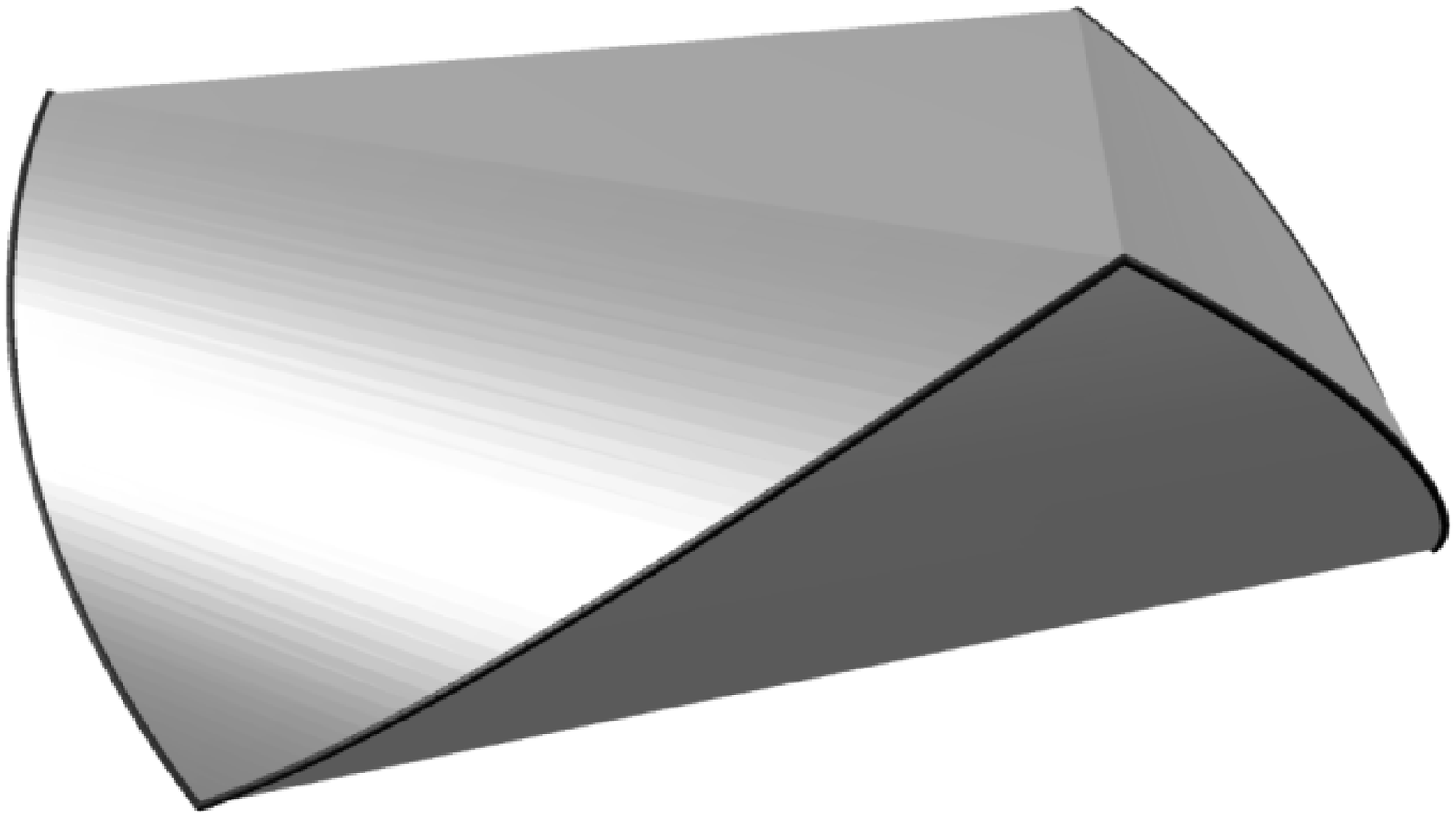}
\caption{Convex surfaces glued from two Euclidean pieces
identified along the boundary. {\em Left:} Disc and equilateral
triangle. {\em Right:} Two Reuleaux triangles (triangles of
constant width); the vertices of one are identified with the
midpoints of the sides of the other.}
\end{figure}

The proof of Alexandrov's theorem presented in this paper provides
a numerical algorithm that constructs a convex polytope from a
given development. This algorithm was implemented by Stefan
Sechelmann, the program is available at {\tt
http://www.math.tu-berlin.de/geometrie/ps/}. Examples of surfaces
constructed with the help of this program are shown in Figure
\ref{fig:Examples}. Two convex polygons $A$ and $B$ of the same
perimeter identified isometrically along the boundary produce a
convex polyhedral metric on the sphere. In the limit, $A$ and $B$
can be any two convex plane figures. The resulting convex body in
$\R^3$ is the convex hull of a space curve that splits its
boundary into two regions with Euclidean metric. This construction
was communicated to us by Johannes Wallner.

There exist numerous generalizations of Alexandrov's theorem: to
surfaces of arbitrary genus, to hyperbolic and spherical
polyhedral metrics, to singularities of negative curvature (see
\cite{Fil} and references therein). Perhaps our approach can be
generalized to these cases.

Note that in our proof we use a special deformation of the metric
implicitely described by the evolution (\ref{eqn:Rule}) of the
curvature. It would be interesting to find a smooth analog of the
flow (\ref{eqn:Rule}) and clarify whether it can be applied to
Weyl's problem.

\section{Generalized convex polytopes and weighted Delaunay triangulations}\label{sec:CGPandDel}
In this section we introduce the notion of a generalized convex
polytope. A generalized convex polytope is a ball equipped with a
Euclidean polyhedral metric that has a convex boundary and a
unique interior vertex.
%In Section \ref{sec:MixVol} we study total curvature of generalized convex polytopes which plays a key role in the proof.
For a fixed Euclidean polyhedral metric $M$ on the sphere, we
consider the space $\P(M)$ of all generalized convex polytopes
with boundary $M$. Our goal here is to introduce coordinates on
$\P(M)$ and to describe its boundary. We also discover a
connection between generalized convex polytopes and weighted
Delaunay triangulations of a sphere with a Euclidean polyhedral
metric. This leads to an explicit description of the space of
weighted Delaunay triangulations.

\subsection{Geodesic triangulations}
Let $M$ be a closed surface with a Euclidean polyhedral metric
with a non-empty singular set $\Sigma$. A \emph{geodesic
triangulation} $T$ of $M$ is a decomposition of $M$ into Euclidean
triangles by geodesics with endpoints in $\Sigma$. By $\V(T)$,
$\E(T)$, $\F(T)$ we denote the sets of vertices, edges, and faces
of $T$, respectively. By definition, $\V(T) = \Sigma$. Note that
in $\E(T)$ there may be multiple edges as well as loops. Triangles
of $\F(T)$ may have self-identifications on the boundary.

Singularities of $M$ are denoted by letters $i,j,k,\ldots$. An
edge of $T$ with endpoints $i$ and $j$ is denoted $ij$, a triangle
with vertices $i$, $j$, $k$ is denoted $ijk$. Note that different
edges or triangles can have the same name and letters in a name
may repeat.

It is not hard to show that for any $M$ there is a geodesic
triangulation $T$. In the case when $M$ is a sphere with all
singularities of positive curvature, this is proved in
\cite[Section 4.1, Lemma 2]{Al05}. Below, examples of geodesic
triangulations are given.

Take two copies of a euclidean triangle and identify them along
the boundary. The space obtained is a geodesically triangulated
sphere with a Euclidean polyhedral metric. There are other
geodesic triangulations of the same surface. Namely, remove from
the triangulation a side of the triangle that forms acute angles
with two other sides and draw instead a geodesic loop based at the
opposite vertex (the loop is formed by two copies of the
corresponding altitude of the triangle).

Let $M$ be the boundary of a cube. Take two opposite faces and
subdivide each of them by a diagonal. The remaining four faces
form a cylinder, both ends of which contain four points from
$\Sigma$. There are infinitely many geodesic triangulations of the
cylinder with vertices in $\Sigma$.

\subsection{Generalized convex polytopes} \label{subsec:GCP}
Let $T$ be a geodesic triangulation of $M$. Let $r =
(r_1,\ldots,r_n)$ be a collection of positive numbers, called
radii, subject to the condition that for every triangle $ijk \in
\F(T)$ there exists a pyramid with base $ijk$ and lengths $r_i,
r_j, r_k$ of side edges. The generalized polytope defined by the
couple $(T, r)$ is a polyhedral complex glued from these pyramids
following the gluing rules in the triangulation. As a consequence,
the apices of all pyramids are glued to one point which we call
the apex of the generalized polytope and denote by $a$.

In a generalized polytope consider the interior edge that joins
the apex to the vertex $i$. Denote by $\omega_i$ the sum of the
dihedral angles of pyramids at this edge. If for some $i$ we have
$\omega_i \ne 2\pi$, then the generalized polytope cannot be
embedded into $\R^3$ and has to be considered as an abstract
polyhedral complex. The quantity
$$
\kappa_i = 2\pi - \omega_i
$$
is called the \emph{curvature} at the corresponding edge. Further,
for every boundary edge $ij \in \E(T)$ of the generalized polytope
denote by $\theta_{ij}$ the sum of the dihedral angles of pyramids
at this edge.

\begin{dfn}
A generalized polytope with boundary $M$ is an equivalence class
of couples $(T,r)$ as above, where two couples are equivalent if
and only if there is an isometry between the resulting polyhedra
which is identical on the boundary and maps the apex to the apex.
A generalized convex polytope is a generalized polytope such that
$\theta_{ij} \le \pi$ for any edge $ij \in \E(T)$ in some (and
hence in any) associated triangulation $T$.
\end{dfn}
For brevity, we sometimes use the word polytope when talking about
generalized convex polytopes. Here are two important examples.

It is not hard to show that a generalized convex polytope with
$\kappa_i = 0$ for all $i$ embeds into $\R^3$ and thus produces a
classical convex polytope. Conversely, take a convex polytope in
$\R^3$, triangulate its non-triangular faces without adding new
vertices, and choose a point in its interior. The triangulation of
the boundary and the distances from the chosen point to the
vertices provide us with a couple $(T,r)$. Another couple $(T',r)$
representing the same generalized convex polytope can be obtained
by choosing other triangulations of non-triangular faces. In
contrast, a different choice of an interior point as the apex
leads to a different generalized polytope.

For any triangulation $T$ take an $R > 0$ and put $r_i = R$ for
every $i \in \V(T)$. If $R$ is large enough, then over every
triangle $ijk$ there exists an isosceles pyramid with all side
edges of length $R$. Thus $(T, r)$ defines a generalized polytope.
It is convex if and only if $T$ is a Delaunay triangulation of
$M$, see Definition \ref{dfn:WDel} or \cite{BS}. This is a special
case of Theorem \ref{thm:Del=Pol} or can be shown directly by a
geometric argument similar to the proof of Lemma
\ref{lem:qTr_prop}.

In our proof of Alexandrov's theorem we take a polytope from the
second example and deform it by changing the radii and keeping the
boundary $M$ fixed to a polytope from the first example. The
explanation why this is possible is rather involved and is based
on the analysis of some unexpected connections to the theory of
mixed volumes.

\subsection{A generalized convex polytope is determined by its radii} \label{subsec:Rad2Triang}
We start our study of the space
\begin{equation} \label{eqn:PS}
\P(M) := \{\mbox{generalized convex polytopes with boundary }M\}.
\end{equation}
We have
$$
\P(M) = \bigcup_T \P^T(M),
$$
where $\P^T(M)$ consists of those generalized convex polytopes
that have a representative of the form $(T,r)$ for the given $T$.
The map $(T,r) \mapsto r$ defines an embedding of $\P^T(M)$ into
$\R^n$, where $n = |\Sigma|$. In this subsection we show that this
map is injective on $\P(M)$ that is, for a given assignment of
radii there is at most one generalized convex polytope.

\begin{prp} \label{prp:Rad2Pol}
Suppose that $P, P' \in \P(M)$ are two generalized convex
polytopes represented by the couples $(T,r)$ and $(T',r')$,
respectively. Then $r=r'$ implies $P=P'$.
\end{prp}
We introduce some technical notions that are used in the proof of
this proposition.

Let us call a \emph{$Q$-function} any function of the form
$$
x \mapsto \left\|x-a\right\|^2 + b
$$
on a subset of a Euclidean space $\R^m$. Here $a \in \R^m$ and $b
\in \R$. A \emph{$PQ$-function} on $M$ is a continuous function
which is a $Q$-function on each triangle of some geodesic
triangulation. To define a $Q$-function on a triangle of $T$ one
develops this triangle onto the plane. By Pythagoras theorem, the
restriction of a $Q$-function to a Euclidean subspace is a
$Q$-function on this subspace. Similarly, the restriction of a
$PQ$-function on $M$ to a geodesic arc is a $PQ$-function of the
arc length. We call a $PQ$-function \emph{$Q$-concave} iff on each
geodesic arc it becomes concave after subtracting the squared arc
length. The difference of any two $Q$-functions is a linear
function, so after subtracting the squared arc length we always
obtain a piecewise linear function on a geodesic.

A $PQ$-function can be spelled out as piecewise quadratic
function, and $Q$-concave as quasiconcave.

\begin{dfn}
Let $T$ be a geodesic triangulation of $M$, and let $r$ be an
assignment of positive numbers to the vertices of $T$. Define the
function
$$
q_{T,r}: M \to \R
$$
as the $PQ$-function that takes value $r_i^2$ at the $i$-th
singularity of $M$, and that is a $Q$-function on every triangle
of $T$.
\end{dfn}

\begin{lem} \label{lem:qTr_prop}
A couple $(T,r)$ represents a generalized convex polytope if and
only if $q_{T,r}$ is a positive $Q$-concave $PQ$-function.
\end{lem}
\begin{proof}. Let $(T,r)$ represent a generalized convex polytope. It is immediate that the function $q_{T,r}$ is the squared distance from the apex. Therefore it is a positive $PQ$-function. To show that it is $Q$-concave, it suffices to consider a short geodesic arc crossing an edge $ij \in \E(T)$. Take the adjacent triangles $ijk$ and $ijl$ and develop them onto the plane. We get a quadrilateral $ikjl$. Choose a point $a \in \R^3$ such that $\left\|a - i\right\| = r_i,\, \left\|a - j\right\| = r_j$, and $\left\|a - k\right\| = r_k$. We claim that $\left\|a - l\right\| \ge r_l$. To see this, compare the pyramid over $ikjl$ with the union of two corresponding pyramids in the polytope, and use the fact that $r_l$ is a monotone increasing function of $\theta_{ij}$, for constant $r_i,\, r_j,\, r_k$. Consider the function $q_{T,r}(x) - \left\|x - a\right\|^2$ on the quadrilateral. It is a piecewise linear function that is identically $0$ on the triangle $ijk$ and non-positive at $l$. Therefore it is concave, which implies $Q$-concavity of $q_{T,r}$ on any short geodesic arc across $ij$. The left to right implication is proved. Note that any of the vertices $i, j, k, l$ may coincide, and even triangles $ijk$ and $ijl$ may be two copies of the same triangle.

Now let $(T,r)$ be any couple such that $q_{T,r}$ is positive and
$Q$-concave. Put $r_i = \sqrt{q_{T,r}(i)}$ for every $i$. Let us
show that for every triangle $ijk \in \F(T)$ there exists a
pyramid over it with side edge lengths $r_i,\, r_j,\, r_k$.
Develop the triangle $ijk$ and extend the function $q_{T,r}$ from
it to a $Q$-function $\ext_{ijk}$ on the whole plane. The desired
pyramid exists if and only if $\ext_{ijk}$ is positive everywhere
on $\R^2$: the minimum value of the function is the square of the
pyramids altitude, the point of minimum is the projection of the
apex. So assume that there is a point $y \in \R^2$ such that
$\ext_{ijk}(y) \le 0$. Join $y$ with some point $x$ inside the
triangle $ijk$ by a straight segment. Now take the corresponding
geodesic arc on $M$, which is possible for $x$ and $y$ in general
position. Inside the triangle $ijk$ we have $q_{T,r} =
\ext_{ijk}$. It follows from the $Q$-concavity of $q_{T,r}$ that
$q_{T,r}(y) \le 0$. This contradiction shows that there is a
generalized polytope corresponding to the couple $(T,r)$. The
proof that it is convex literally follows the argument from the
first part of the proof. \end{proof}

\begin{proof}{\it of Proposition \ref{prp:Rad2Pol}.}
Let us show that a generalized convex polytope is determined by
the function $q_{T,r}$:
\begin{equation} \label{eqn:qT}
(T,r) \sim (T',r') \Leftrightarrow q_{T,r} = q_{T',r'}.
\end{equation}
The implication from the left to the right is obvious. To prove
the inverse implication, consider a common geodesic subdivision
$T''$ of $T$ and $T'$, where $T''$ is allowed to have vertices
outside $\Sigma$. This yields decompositions of the both
generalized polytopes $(T,r)$ and $(T',r')$ into smaller pyramids.
The equality $q_{T,r} = q_{T',r'}$ implies that the pyramids over
same triangle of $T''$ in both decompositions are equal.

Now it suffices to show that if for some $r$ there are two
triangulations $T$ and $T'$ such that both functions $q_{T,r}$ and
$q_{T',r}$ are $Q$-concave, then $q_{T,r} = q_{T',r}$.

Take any point $x \in M$. Let $ijk \in \F(T)$ be a triangle that
contains $x$. Since $q_{T,r}$ is a $Q$-function on $ijk$, the
function $q_{T',r} - q_{T,r}$ is piecewise linear and concave on
$ijk$. Since it vanishes at the vertices, it is non-negative
everywhere in $ijk$. Thus we have $q_{T',r}(x) \ge q_{T,r}(x)$ for
any $x \in M$. Similarly, $q_{T,r}(x) \ge q_{T',r}(x)$. Therefore
the functions are equal.
\end{proof}

\subsection{The space of generalized convex polytopes} \label{subsec:CGP}
By Proposition \ref{prp:Rad2Pol}, radii can be considered as
coordinates on the space $\P(M)$ which thus becomes a subset of
$\R^n$. It is more convenient to take squares of radii as
coordinates; so in the future we identify $\P(M)$ with its image
under the embedding
\begin{eqnarray}
\P(M) & \to & \R^n, \label{eqn:Pemb}\\
(T,r) & \mapsto & (r^2_1,\ldots,r^2_n). \nonumber
\end{eqnarray}
In this subsection we derive a system of inequalities which
describes $\P(M)$.

By Lemma \ref{lem:qTr_prop}, we have
$$
\P(M) = \{\mbox{positive }Q\mbox{-concave } PQ\mbox{-functions on
} M\}.
$$
First, we look at a larger space:
$$
\D(M) = \{Q\mbox{-concave } PQ\mbox{-functions on } M\}.
$$
For any map $q: \Sigma \to \R$ and any geodesic triangulation $T$
of $M$ denote by
$$
\widetilde{q_T}: M \to \R
$$
the $PQ$-extension of $q$ with respect to $T$. Note that
$\widetilde{q_T}$ relates to $q$ as $q_{T,r}$ relates to $r^2$.
The space $\D(M)$ can be identified with those $q = (q_i)_{i \in
\Sigma}$ for which there exists $T$ such that $\widetilde{q_T}$ is
$Q$-concave. By $q_i$ we denote the value $q(i)$.

A Euclidean quadrilateral $ikjl$ in $M$ is a region bounded by
simple geodesic arcs $ik, kj, jl, li$ and without singularities in
the interior. Any two of the boundary arcs can coincide but cannot
intersect. A Euclidean quadrilateral can be developed onto the
plane, after resolving possible self-identifications on the
boundary. A Euclidean triangle $iji$ arises when there is a closed
geodesic based at $i$ that encloses a unique singularity $j$.

Define the function $\ext_{ikl}$ on the quadrilateral $ikjl$ as
the $Q$-function that takes values $q_i,\, q_k,\, q_l$ at the
respective vertices.

Now we describe the \emph{flip algorithm}. Let $ij$ be an edge in
a geodesic triangulation $T$. If $ij$ is adjacent to two different
triangles $ijk$ and $ijl$, then we have a Euclidean quadrilateral
$ikjl$. If $ikjl$ is strictly convex, then we can transform the
triangulation $T$ by replacing the diagonal $ij$ through the
diagonal $kl$ in $ikjl$. This transformation is called a
\emph{flip}. The flip algorithm in a general setting works as
follows. Assume that we have a rule to say for any edge $ij \in
\E(T)$ whether $ij$ is ``bad'' or ``good''. The goal is to find a
triangulation where all edges are good. Start with an arbitrary
triangulation. If there is a bad edge that forms a diagonal of a
strictly convex quadrilateral, then flip it. Look for a bad edge
in the new triangulation, flip it if possible and so on. If any
bad edge can be flipped and flipping bad edges cannot continue
infinitely, then the flip algorithm yields a triangulation with
only good edges.

\begin{prp} \label{prp:DelSpace}
The set $\D(M)$ is a convex polyhedron which is the solution set
of a system of linear inequalities of the form:
\begin{eqnarray}
q_j & \ge & \ext_{ikl}(j), \label{eqn:ConcQuadr}\\
q_j & \ge & q_i - \l^2_{ij}. \label{eqn:CloGeod}
\end{eqnarray}
There is one equation of the form (\ref{eqn:ConcQuadr}) for each
Euclidean quadrilateral $ikjl$ with the angle at $j$ greater or
equal $\pi$, and one equation of the form (\ref{eqn:CloGeod}) for
each Euclidean triangle $iji$. By $\l_{ij}$ we denote the length
of the edge $ij$.
\end{prp}

\begin{figure}[ht]
\centerline{\input{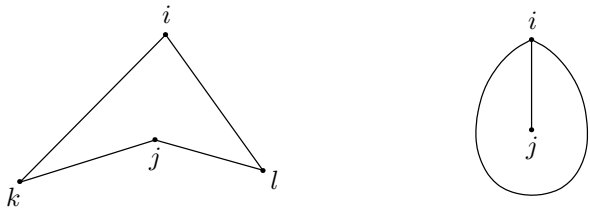}} \caption{Every
non-strongly convex quadrilateral $ikjl$ gives rise to an
inequality (\ref{eqn:ConcQuadr}), every triangle $iji$ gives rise
to an inequality (\ref{eqn:CloGeod}).}
\end{figure}

\begin{proof}. Let us show that the right hand side of (\ref{eqn:ConcQuadr}) is a linear function of $q_i,\, q_k,\, q_l$. Develop $ikjl$ onto the plane. When we vary $q_i$, the function $\ext_{ikl}$ changes by a linear function which vanishes along the line $kl$. Thus $\ext_{ikl}(j)$ depends linearly on $q_i$, as well as on $q_k$ and $q_l$. Moreover, the inequality (\ref{eqn:ConcQuadr}) has the form
$$
q_i \ge \lambda q_j + \mu q_k + \nu q_l + c,
$$
where $\lambda + \mu + \nu = 1$, $\lambda,\, \mu,\, \nu \ge 0$ are
barycentric coordinates of $j$ with respect to $i,\, k$, and $l$.

Let us show that conditions (\ref{eqn:ConcQuadr}) and
(\ref{eqn:CloGeod}) are fulfilled for any assignment $q = (q_i)_{i
\in \Sigma}$ that posesses a $Q$-concave extension
$\widetilde{q}$. For (\ref{eqn:ConcQuadr}), consider the function
$\widetilde{q} - \ext_{ikl}$ on the quadrilateral $ikjl$. It is
concave and vanishes at $i,\, k$, and $l$. Hence it is
non-negative on the triangle $ikl$, in particular $q_j -
\ext_{ikl}(j) \ge 0$. For (\ref{eqn:CloGeod}), consider the
function $f(x) = q(x) - \left\|x - j\right\|^2$ on the triangle
$iji$. It is geodesically concave and takes values $q_j$ at $j$
and $q_i - \l^2_{ij}$ at $i$. Hence $f(x) \ge q_i - \l^2_{ij}$ for
any $x$ on the edge $ii$. Due to concavity, $f$ is also greater or
equal to $q_i - \l^2_{ij}$ at any regular point of the triangle
$iji$. By continuity, $f(j) \ge q_i - \l^2_{ij}$.

Now we prove sufficiency of the conditions (\ref{eqn:ConcQuadr})
and (\ref{eqn:CloGeod}). Let $T$ be any geodesic triangulation.
Call an edge of $T$ \emph{good}, if the function $\widetilde{q_T}$
is $Q$-concave across this edge, otherwise call an edge
\emph{bad}. Our aim is to find a triangulation $T$ with only good
edges. Let us show that any bad edge can be flipped. Indeed, if
$ij$ is adjacent to only one triangle, then the inequality
(\ref{eqn:CloGeod}) implies that $ij$ is good. If it is adjacent
to two triangles that form a concave or non-strictly convex
quadrilateral, then the inequality (\ref{eqn:ConcQuadr}) implies
that $ij$ is good. Thus any bad edge is a diagonal of a strictly
convex quadrilateral and therefore can be flipped.

Let us show that the flip algorithm terminates. Note that when a
bad edge is flipped, the function $\widetilde{q_T}$ increases
pointwise. If there are only finitely many geodesic triangulations
of $M$, this suffices. If not, we use the fact that for any $C$
there are only finitely many triangulations of $M$ with all
geodesics of length less than $C$. Thus it suffices to show that
no long edges can appear when performing the algorithm. Let $T'$
be a triangulation with an edge $ij$. Function
$\widetilde{q_{T'}}$ restricted to $ij$ has the form $x^2 + ax +
b$ and takes values $q_i$ and $q_j$ at the endpoints. There exists
a constant $C$ such that if the length of $ij$ is greater than
$C$, the minimum value of $\widetilde{q_{T'}}$ on $ij$ is smaller
than $\min_{x \in M} \widetilde{q_T}(x)$. This means that no edge
of length greater than $C$ can appear between the vertices $i$ and
$j$ if we start the flip algorithm from the triangulation $T$.
\end{proof}

Let $\D^T(M)$ be the set of those functions $q \in \D(M)$ which
are $PQ$ with respect to the triangulation $T$. Then we have a
decomposition
\begin{equation} \label{eqn:DecDM}
\D(M) = \cup_T \D^T(M).
\end{equation}

\begin{prp} \label{prp:PropDel}
For any surface $M$ with a polyhedral metric the polyhedron
$\D(M)$ of $Q$-concave $PQ$-functions on $M$ has the following
properties.
\begin{enumerate}
\item \label{it:1} The origin $0 \in \R^n$ is an interior point of
$\D(M)$. \item \label{it:2} Put $\D_0(M) = \{q \in \D(M)|\: \sum_i
q_i = 0\}$. Then the polyhedron $\D(M)$ is the direct sum
\begin{equation} \label{eqn:Decomp}
\D(M) = \D_0(M) \oplus L,
\end{equation}
where $L$ is the one-dimensional subspace of $\R^n$ spanned by the
vector $(1,1,\ldots,1)$. Every $\D^T(M)$ is decomposed in the same
way: if $\D_0^T(M) = \{q \in \D^T(M)|\: \sum_i q_i = 0\}$, then
$\D^T(M) = \D_0^T(M) \oplus L$. \item \label{it:3} Every $\D^T(M)$
is a convex polyhedron. \item \label{it:4} The decomposition
(\ref{eqn:DecDM}) is locally finite. \item \label{it:5} For $q \in
\D(M)$, let $\widetilde{q}: M \to \R$ denote the $Q$-concave
$PQ$-extension of $q$. Then $\widetilde{q}$ depends continuously
on $q$ in the $L^\infty$-metric.
\end{enumerate}
\end{prp}
\begin{proof}. If we put $q_i = 0$ for every $i$, then all of the inequalities (\ref{eqn:ConcQuadr}) and (\ref{eqn:CloGeod}) are fulfilled and are strict. This implies property \ref{it:1}.

Fix a triangulation $T$. If $q$ is a $PQ$-function with respect to
$T$, then it is $Q$-concave across an edge $ij$ if and only if
\begin{equation} \label{eqn:DT}
q_l \le \ext_{ijk}(l),
\end{equation}
where $ijk$ and $ijl$ are triangles of $T$. The inequality
(\ref{eqn:DT}) remains valid if we change all of the $q_i$'s by
the same amount. Thus $D^T(M)$, and also the whole $D(M)$, is
invariant under translations along $L$. This proves
property~\ref{it:2}.

The set $D^T(M)$ is a convex polyhedron because it is the set of
solutions of a linear system (\ref{eqn:DT}). This shows
\ref{it:3}.

Property \ref{it:4}. means that any bounded subset $U$ in $\R^n$
has a non-empty intersection with only finitely many of the
$\D^T(M)$'s. For this, it suffices to show that the edge lenghts
in all of the triangulations $T$ such that $\D^{T}(M) \cap U \ne
\emptyset$ are uniformly bounded from above. Fix any triangulation
$T_0$. For any $q \in U$ consider the functions
$\widetilde{q_{T_0}}$ and $\widetilde{q_T}$, where $T$ is such
that $\widetilde{q_T}$ is $Q$-concave. We have
$\widetilde{q_{T_0}} \le \widetilde{q_T}$. Since $U$ is bounded,
functions $\widetilde{q_{T_0}}$ are uniformly bounded from below.
Take any edge in $T$. The values of $q$ at the endpoints of the
edge are less or equal $\max_{q \in U,\, i \in \Sigma} q_i$. This,
together with the uniform lower bound on $q$, implies that the
edge cannot be too long.

In every $\D^T(M)$, $\widetilde{q}$ depends continuously on $q$.
The continuity of $q \mapsto \widetilde{q}$ on the whole $\D(M)$
follows from property \ref{it:4}. \end{proof}

\begin{prp} \label{prp:P(M)}
The space $\P(M)$ of convex generalized polytopes with boundary
$M$ has the following properties.
\begin{enumerate}
\item \label{itt:1} The point $(C,C,\ldots,C)$ is an interior
point of $\P(M)$ for a sufficiently large $C$. \item \label{itt:3}
If $q \in \P(M)$ is an interior point of $\D(M)$, then $q$ is an
interior point of $\P(M)$ as well. \item \label{itt:4} $ \P^T(M) =
\D^T(M) \cap \{q_i > 0 \, \forall i \in \V(T) \} \cap \{CM_{ijk} >
0 \, \forall ijk \in \F(T)\}, $ where
$$
CM_{ijk} = \left|
\begin{array}{ccccc}
0 & 1 & 1 & 1 & 1\\
1 & 0 & q_i & q_j & q_k \\
1 & q_i & 0 & \l_{ij}^2 & \l_{ik}^2 \\
1 & q_j & \l_{ij}^2 & 0 & \l_{jk}^2 \\
1 & q_k & \l_{ik}^2 & \l_{jk}^2 & 0
\end{array}
\right|
$$
is the Cayley-Menger determinant.
\end{enumerate}
\end{prp}
\begin{proof}. Recall that
$$
\P(M) = \{q \in \D(M)|\: \widetilde{q} > 0\}.
$$
By Proposition \ref{prp:PropDel}, the point $(q_i = C)$ belongs to
$\D(M)$ for any $C$. Since adding a constant to all of $q_i$
results in adding a constant function to $\widetilde{q}$, for a
large $C$ the $PQ$-extension of $(q_i = C)$ is everywhere
positive. That is, $(q_i = C) \in \P(M)$.

By property \ref{it:5}. of Proposition \ref{prp:PropDel}, the
function $\widetilde{q}$ depends continuously on $q$. Therefore
the condition $\widetilde{q} > 0$ is an open condition on $q$.
This implies property \ref{itt:3}. Also it completes the proof of
property \ref{itt:1}, since ${q_i = C}$ is an interior point of
$\D(M)$.

Positivity of the Cayley-Menger determinant $CM_{ijk}$ is
equivalent to the existence of a pyramid over the triangle $ijk$
with side edges of lengths $\sqrt{q_i}, \sqrt{q_j}, \sqrt{q_k}$.
This implies property \ref{itt:4}. \end{proof}

Note that $q$ is a boundary point of $\D(M)$ if and only if one of
the inequalities (\ref{eqn:ConcQuadr}) and (\ref{eqn:CloGeod})
becomes an equality. Geometrically this means that a generalized
convex polytope $P$ is on the boundary of $\P(M)$ if and only if
$P$ has a non-strongly convex face or an isolated vertex. Here by
faces of $P$ we mean connected regions in $M$ after erasing flat
edges of $P$.

\subsection{Weighted Delaunay triangulations} \label{subsec:WDelaunay}
Here we explain how generalized convex polytopes are related to
weighted Delaunay triangulations of polyhedral surfaces.

\begin{dfn} \label{dfn:WDel}
Let $M$ be a polyhedral surface with singular set $\Sigma$. Let $V
\supset \Sigma$ be a finite non-empty subset that has at least one
point in every boundary component of $M$. For a function $q: V \to
\R$ and a geodesic triangulation $T$ with $\V(T) = V$ denote by
$\widetilde{q_T}$ the $PQ$-extension of $q$ to $M$ with respect to
$T$. A couple $(T,q)$ is called a \emph{weighted Delaunay
triangulation} of $(M, V)$ with weights $q$ if $\widetilde{q_T}$
is $Q$-concave.

A Delaunay triangulation is a weighted Delaunay triangulation with
all weights equal.
\end{dfn}

\begin{dfn}
A weighted Delaunay tesselation $(\overline T, q)$ is obtained
from a weighted Delaunay triangulation $(T,q)$ by removing from
$T$ all inessential edges. An edge $ij \in \E(T)$ is called
inessential if the corresponding function $\widetilde{q_T}$ is a
$Q$-function in the neighborhood of any interior point of $ij$.
\end{dfn}

Originally weighted Delaunay triangulations were defined for $V
\subset \R^2$, and $M$ the convex hull of $V$, see \cite{AK},
\cite{Ed}, \cite{For}. It can be shown that in this case
Definition \ref{dfn:WDel} is equivalent to the classical one.
Weighted Delaunay triangulations of polyhedral surfaces appear in
\cite{Gli}, \cite{Sch2}. However, Definition \ref{dfn:WDel} seems
to be new.

\begin{thm} \label{thm:Del=Pol}
If $M$ is the sphere and $V = \Sigma$, then weighted Delaunay
tesselations of $(M,V)$ with positive associated functions
$\widetilde{q_T}$ are in one-to-one correspondence with
generalized convex polytopes with boundary $M$.
\end{thm}
\begin{proof}. This is essentially a reformulation of Lemma \ref{lem:qTr_prop}. By (\ref{eqn:qT}), generalized convex polytopes with boundary $M$ are in one-to-one correspondence with positive $Q$-concave $PQ$-functions on $M$. To any such function $f$ there canonically corresponds a tesselation of $M$ into regions where $f$ is quadratic. Since the function $f$ is $Q$-concave, the tesselation together with values of $f$ at $\Sigma$ as weights is a weighted Delaunay tesselation.
\end{proof}

Note that a weighted Delaunay triangulation subsumes an assignment
of weights to the vertices. If $T$ is a triangulation for which
there exist weights $q$ such that $(T,q)$ is a weighted Delaunay
triangulation, then $T$ is sometimes called regular or coherent
triangulation.

\begin{dfn}
Let $M$ and $V$ be as in Definition \ref{dfn:WDel}. Define the
space of \emph{admissible weights} $\D(M,V)$ as the set of maps
$q: V \to \R$ such that there exists a weighted Delaunay
triangulation of $M$ with weights $q$.
\end{dfn}

\begin{thm} \label{thm:Delaunay}

\begin{enumerate}
\item For any $(M,V)$ the space of admissible weights $\D(M,V)$ is
a convex polyhedron in $\R^V$ defined by inequalities
(\ref{eqn:ConcQuadr}) and (\ref{eqn:CloGeod}). \item For any
weight $q \in \D(M,V)$ the weighted Delaunay tesselation with
weights $q$ is unique. \item A weighted Delaunay triangulation can
be found via the flip algorithm starting from any triangulation of
$(M,V)$. \item Properties from Proposition \ref{prp:PropDel} hold.
\end{enumerate}
\end{thm}

\begin{proof}. Statement 1. generalizes Proposition \ref{prp:DelSpace}. In Proposition \ref{prp:DelSpace} we assumed that $M$ is a closed surface and $V = \Sigma$. But it is easy to see that the proof works equally well without these assumptions. Statement 3. is contained in the proof of Proposition \ref{prp:DelSpace}.

Statement 2. is essentially Proposition \ref{prp:Rad2Pol}.

Finally, for statement 4. nothing changes in the proof of
Proposition~\ref{prp:PropDel}.\!\!\!\!\!~\end{proof}

The weighted Delaunay tesselation with all weights equal is called
the {\em Delaunay tesselation}. Theorem \ref{thm:Delaunay} implies
that for any $(M,V)$ the Delaunay tesselation is unique and a
Delaunay triangulation can be obtained via the flip algorithm.

Delaunay triangulations of polyhedral surfaces were used as a
technical tool in the study of moduli spaces of Euclidean
polyhedral metrics on surfaces, \cite{Bow}, \cite{Riv}. Recently
they were applied to define the discrete Laplace-Beltrami operator
on polyhedral surfaces \cite{BS}.

\section{Total scalar curvature and volume of the dual} \label{sec:MixVol}
In Section \ref{sec:CGPandDel} we have shown that the radii
$(r_i)_{i \in \Sigma}$ can serve as coordinates on the space
$\P(M)$ of generalized convex polytopes with boundary $M$. Hence
the curvatures $(\kappa_i)_{i \in \Sigma}$ become functions of the
radii. In our proof of Alexandrov's theorem we will deform the
radii so that at the end of the deformation all of the curvatures
vanish. The goal of this Section is to prove non-degeneracy of the
Jacobian $\left(\frac{\partial \kappa_i}{\partial r_j}\right)$
under certain restrictions on the polytope, Corollary
\ref{cor:NonDeg}. By the inverse function theorem, this allows to
realize any infinitesimal deformation of curvatures by an
infinitesimal deformation of radii.

\subsection{Total scalar curvature of a generalized polytope}
\begin{dfn}
Let $P = (T,r)$ be a convex generalized polytope. The total scalar
curvature of $P$ is
$$
H(P) = \sum_{i \in \V(T)} r_i \kappa_i + \sum_{ij \in \E(T)}
\l_{ij}(\pi - \theta_{ij}),
$$
where $\kappa_i = 2\pi - \omega_i$ with $\omega_i$ the sum of the
dihedral angles of pyramids at the $i$-th radial edge, $\l_{ij}$
is the length of the edge $ij \in \E(T)$, and $\theta_{ij}$ is the
dihedral angle of the polytope $P$ at the edge $ij$.
\end{dfn}
If triangulation $T$ is not unique, the value $H(P)$ does not
depend on it, because $\pi - \theta_{ij} = 0$ for a flat edge
$ij$.
\begin{prp} \label{prp:CurvDer}
For a fixed Euclidean polyhedral metric $M$, function $H$ is of
class $C^2$ on $\P(M)$. Its partial derivatives are:
%$$
%\frac{\partial \kappa_i}{\partial r_j} = \frac{\partial^2 H}{\partial r_i \partial r_j} = \frac{\partial \kappa_j}{\partial r_i}
%$$
\begin{eqnarray}
\frac{\partial H}{\partial r_i} & = & \kappa_i, \label{eqn:DmDr}\\
%\frac{\partial \kappa_i}{\partial r_j} & = & 0 \mbox{ if } ij \notin T \\
\frac{\partial \kappa_i}{\partial r_j} & = & \frac{\cot \alpha_{ij} + \cot \alpha_{ji}}{\l_{ij} \sin \rho_{ij} \sin \rho_{ji}}, \label{eqn:KiRj}\\%\mbox{ if } ij \in T\\
\frac{\partial \kappa_i}{\partial r_i} & = & - \sum_{j \ne i} \cos
\phi_{ij} \frac{\partial \kappa_i}{\partial r_j}. \label{eqn:KiRi}
\end{eqnarray}
Here $\alpha_{ij}$ and $\alpha_{ji}$ are the dihedral angles of
the pyramids at the edge $ij$, thus $\alpha_{ij} + \alpha_{ji} =
\theta_{ij}$; $\rho_{ij}$ is the angle at the vertex $i$ in the
triangle $aij$; $\phi_{ij}$ is the angle at $a$ in the same
triangle. If $ij$ is not an edge of $T$, then $\frac{\partial
\kappa_i}{\partial r_j} = 0$.
\end{prp}

\begin{figure}[ht]
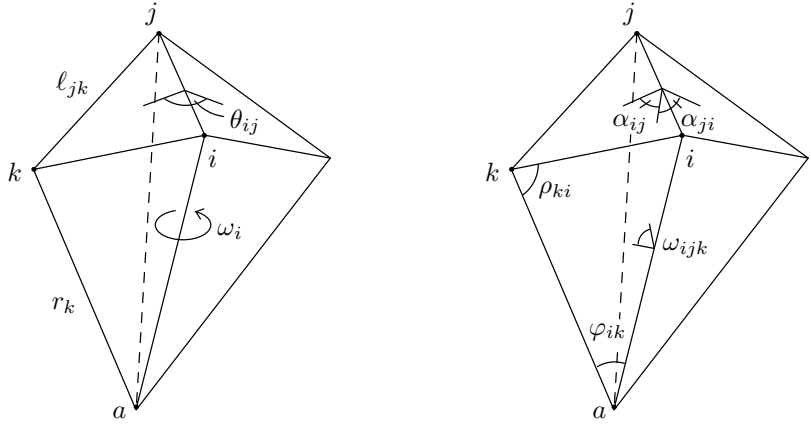

\centerline{\input{Figures/fig02.pstex_t}\hspace{2cm}\input{Figures/fig03.pstex_t}}
\caption{Angles and lengths in a generalized convex polytope.}
\label{fig:fig02}
\end{figure}

\begin{proof}. Let $P$ be strongly convex, that is let $\theta_{ij} < \pi$ for all $ij \in T$. In this case all polytopes in some neighborhood of $P$ have the same triangulation $T$ on the boundary.

For $ijk \in \F(T)$ consider the corresponding pyramid. Denote by
$\omega_{ijk}$ its dihedral angle at the edge $ai$. By
Schl{\"a}fli formula \cite{Mil},
$$
r_i d\omega_{ijk} + r_j d\omega_{jik} + r_k d\omega_{kij} +
\l_{ij} d\alpha_{ij} + \l_{jk} d\alpha_{jk} + \l_{ki} d\alpha_{ki}
= 0.
$$
Summing this up over all $ijk \in \F(T)$ yields
$$
\sum_{i \in \V(T)} r_i d\omega_i + \sum_{ij \in \E(T)} \l_{ij}
d\theta_{ij} = 0
$$
which implies formula (\ref{eqn:DmDr}).

For $i \in \V(T)$ denote
\begin{equation} \label{eqn:LinkDef}
\lk i = \{j \in \V(T)|\: ij \in \E(T)\}.
\end{equation}
The angle $\omega_i$ and its parts $\omega_{ijk}$ can be viewed as
functions of the angles $(\rho_{ij})_{j \in \lk i}$. By Leibniz
rule we have
$$
\begin{array}{rcl}
\frac{\partial \omega_{ijk}}{\partial r_j} & = & \frac{\partial \omega_{ijk}}{\partial \rho_{ij}}\frac{\partial \rho_{ij}}{\partial r_j},\\
\frac{\partial \omega_{ijk}}{\partial r_i} & = & \frac{\partial
\omega_{ijk}}{\partial \rho_{ij}}\frac{\partial
\rho_{ij}}{\partial r_i} + \frac{\partial \omega_{ijk}}{\partial
\rho_{ik}}\frac{\partial \rho_{ik}}{\partial r_i}.
\end{array}
$$
All partial derivatives on the right hand side are derivatives of
the angles with respect to the side lengths in spherical and
Euclidean triangles. Substituting the formulas of Lemma
\ref{lem:SphBasic}, we get the formulas (\ref{eqn:KiRj}) and
(\ref{eqn:KiRi}).

Let $P$ be not strongly convex. Recall that $\P^T(M)$ denotes the
subset of $\P(M)$ formed by the polytopes with an associated
triangulation $T$. By Property \ref{it:3}. from Proposition
\ref{prp:PropDel}, $\P^T(M)$ are separated by affine hyperplanes
in coordinates $(q_i) = (r^2_i)$. By the argument for a strongly
convex $P$, there are directional derivatives of $H$ at $P$ in all
directions. Since $\kappa_i$ is well-defined and continuous on
$\P(M)$, it follows that function $H$ is of class $C^1$. The same
argument applies to functions $\kappa_i$. The right hand sides of
formulas (\ref{eqn:KiRj}) and (\ref{eqn:KiRi}) are well-defined
and continuous because all of the summands that correspond to flat
edges vanish.
\end{proof}

In fact, the formulas of Proposition \ref{prp:CurvDer} are valid
only when the triangulation $T$ contains neither loops nor
multiple edges. Otherwise the formulation is somewhat cumbersome
and looks as follows. For an oriented edge $e \in \E(T)$ denote by
$a(e), b(e) \in \V(T)$ its initial, respectively terminal, vertex.
Then we have
$$
d\kappa_i = \sum_{a(e) = i} \frac{\cot \alpha_e + \cot
\alpha_{-e}}{\sin \rho_e} d\rho_e,
$$
\begin{equation} \label{eqn:Complic}
\frac{\partial\rho_e}{\partial r_j} = \left\{ \begin{array}{ll}
                \frac{1}{\ell_e \sin\rho_e}&\mbox{ for }a(e) \ne j = b(e)\\
                -\frac{\cos\phi_e}{\ell_e\sin\rho_{-e}}&\mbox{ for }a(e) = j\ne b(e)\\
                \frac{\ell_e}{2r^2_i\sin\rho_e}&\mbox{ for }a(e) = j = b(e).
                                              \end{array}\right.
\end{equation}

\subsection{Generalized polyhedra} \label{subsec:GenPolyh}
Convex polytope and convex polyhedron are two basic notions of the
classical theory of polytopes. The former is defined as the convex
hull of a finite number of points, and the latter as the
intersection of a finite number of half-spaces. If we take a
sphere with the center inside a convex polytope $P$, then the
corresponding polar dual $P^*$ is a bounded convex polyhedron, and
vice versa. On the other hand, there is a fundamental theorem
saying that the class of bounded convex polyhedra coincides with
the class of convex polytopes. For more information see
\cite{Zie}.

Here we dualize the notion of generalized convex polytope from
Section \ref{sec:CGPandDel} by introducing generalized convex
polyhedra. Informally speaking, a generalized convex polytope has
curvatures concentrated along the segments joining the apex to the
vertices. A generalized convex polyhedron has curvatures
concentrated along the altitudes drawn from the apex to the faces.
A common extension of these two notions is a polyhedral complex
that has curvatures both along the radial edges and along the
altitudes. This class of complexes is self-dual; however, we have
no need of such generalization.

Generalized convex polyhedra are introduced in Definitions
\ref{dfn:GSF} and \ref{dfn:Polyh}. In order to motivate the
abstract constructions let us look closer at the classical convex
polyhedra.

Let us start by recalling the notion of the normal fan. Let $Q$ be
a bounded 3-dimensional convex polyhedron. To a face $Q_i$ of $Q$
associate the ray in $\R^3$ based at the origin and spanned by the
outer normal to $Q_i$. If two faces $Q_i$ and $Q_j$ have an edge
in common, then denote this edge by $Q_{ij}$ and associate to it
the positive span of the rays associated with $Q_i$ and $Q_j$. The
flat angles thus obtained subdivide $\R^3$ into 3-dimensional
cones which are the normal cones at the vertices of $Q$. The
resulting complex of cones is called the \emph{normal fan} of $Q$.
If exactly three edges meet at every vertex of $Q$, then all
vertex cones are trihedral angles, and the normal fan is called
\emph{simplicial}. Every normal fan can be subdivided to become
simplicial, which corresponds to formally introducing in $Q$
additional edges of zero length. The spherical section of a
simplicial fan is a geodesic triangulation of $\Sph^2$. Clearly,
it completely determines the fan.

\begin{dfn} \label{dfn:GSF}
A generalized simplicial fan is a couple $(T,\phi)$, where $T$ is
a triangulation of the sphere, and
$$
\phi: \E(T) \to (0, \pi)
$$
is a map such that for each triangle $ijk \in \F(T)$ there exists
a spherical triangle with side lengths $\phi_{ij},\, \phi_{jk},\,
\phi_{ik}$.
\end{dfn}

A simplex $OABC$ in $\R^3$ is called an \emph{orthoscheme}, if the
vectors $OA$, $AB$, $BC$ are pairwise orthogonal. Let $Q \subset
\R^3$ be a convex polyhedron. Draw perpendiculars from a point $O
\in \R^3$ to the planes of the faces of $Q$. Let $A_i$ be the foot
of the perpendicular to the face $Q_i$. In the plane of $Q_i$ draw
perpendiculars $A_i A_{ij}$ to the edges $Q_{ij}$. By $A_{ijk}$
denote the vertex of $Q$ common to the faces $Q_i, Q_j$ and $Q_k$.
The collection of orthoschemes $O A_i A_{ij} A_{ijk}$ is called
the \emph{orthoscheme decomposition} of $Q$. Combinatorially the
orthoscheme decomposition is the barycentric decomposition of $Q$.
Note although that the foots of perpendiculars can lie outside the
faces or edges, which makes the geometric picture complicated. Let
$h_i$, $h_{ij}$, $h_{ijk}$ be the signed lengths of the segments
$O A_i$, $A_i A_{ij}$, $A_{ij} A_{ijk}$, respectively. The length
$h_i$ is negative if the point $O$ and the polytope $Q$ are on the
opposite sides of the face $Q_i$; similarly for $h_{ij}$ and
$h_{ijk}$. We call the numbers $h_i$, $h_{ij}$, $h_{ijk}$ the
\emph{parameters of the orthoscheme} $O A_i A_{ij} A_{ijk}$.

\begin{figure}[ht]
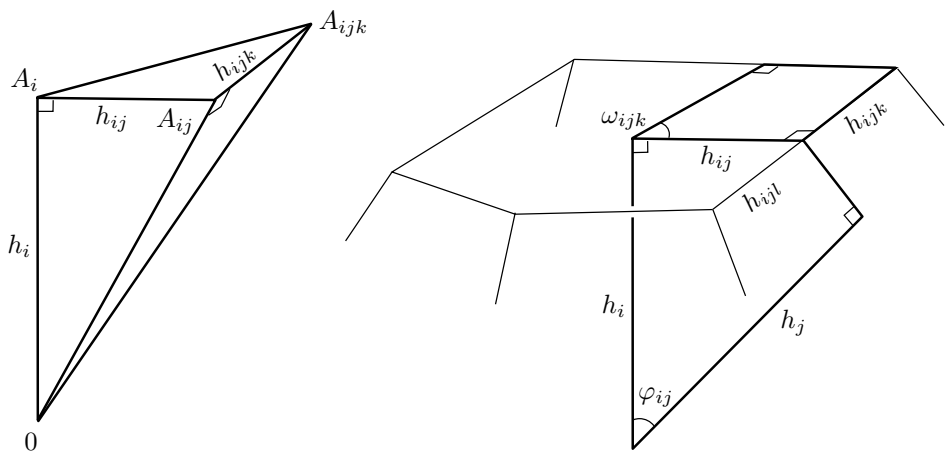

\centerline{\input{Figures/fig04.pstex_t}\hspace{.3cm}\input{Figures/fig05.pstex_t}}
\caption{The orthoscheme decomposition of a convex polyhedron.}
\end{figure}

Clearly, the normal fan and the altitudes $h_i$ completely
determine the polyhedron $Q$. In particular, the parameters
$h_{ij}$ and $h_{ijk}$ can be expressed through $(h_i)$ and the
angles of the fan. One can compute
\begin{equation} \label{eqn:Dhij}
h_{ij} = \frac{h_j - h_i \cos \phi_{ij}}{\sin \phi_{ij}},
\end{equation}
where $\phi_{ij}$ is the angle between the outer normals to faces
$Q_i$ and $Q_j$. Similarly,
\begin{equation} \label{eqn:Dhijk}
h_{ijk} = \frac{h_{ik} - h_{ij} \cos \omega_{ijk}}{\sin
\omega_{ijk}},
\end{equation}
where $\omega_{ijk}$ is the exterior angle of the face $Q_i$ at
the vertex $A_{ijk}$ or, equivalently, the angle at the vertex $i$
of the triangle $ijk$ in the spherical section of the normal fan.
The length $\l^*_{ij}$ of the edge $Q_{ij}$ can be expressed as
\begin{equation} \label{eqn:Dl}
\l^*_{ij} = h_{ijk} + h_{ijl},
\end{equation}
where $Q_k$ and $Q_l$ are the faces that bound the edge $Q_{ij}$.
Therefore the sum $h_{ijk} + h_{ijl}$ has to be non-negative,
although some of the parameters $h_{ijk}$ are allowed to be
negative. Note also that
\begin{equation} \label{eqn:ijk=jik}
h_{ijk} = h_{jik}
\end{equation}
which is not obvious from (\ref{eqn:Dhijk}), but is clear from the
geometric meaning of $h_{ijk}$.

A generalized convex polyhedron is defined as a generalized
simplicial fan together with the set of altitudes $h_i$. More
intuitively, it is a complex of signed orthoschemes with
parameters computed by formulas (\ref{eqn:Dhij}) and
(\ref{eqn:Dhijk}). Moreover, we require the edge lengths
(\ref{eqn:Dl}) to be non-negative.

\begin{dfn} \label{dfn:Polyh}
Let $(T, \phi)$ be a generalized simplicial fan, and $h$ be a map
$\V(T) \to \R$. A triple $Q = (T, \phi, h)$ is called a {\rm
generalized convex polyhedron}, if $\l^*_{ij} \ge 0$ for every
edge $ij \in \E(T)$, where $\l^*_{ij}$ is defined by
(\ref{eqn:Dhij}), (\ref{eqn:Dhijk}), and (\ref{eqn:Dl}). Here
$\omega_{ijk}$ is the angle at the vertex $i$ in the spherical
triangle $ijk$ of the fan $(T, \phi)$.

For every $i \in \V(T)$ define the area of the $i$-th face of $Q$
by
\begin{equation} \label{eqn:DArea}
F_i = \frac{1}{2} \sum_j h_{ij} \l^*_{ij}.
\end{equation}
The volume of the polyhedron $Q$ is defined by
\begin{equation} \label{eqn:DVol}
\vol = \frac{1}{3} \sum_i h_i F_i.
\end{equation}
\end{dfn}
Note that the right hand side of the formula (\ref{eqn:DVol}) is
nothing else than the sum of the signed volumes of the
orthoschemes.

The space $\Q(T, \phi)$ of all generalized convex polyhedra with a
given normal fan $(T,\phi)$ is a convex polyhedron in $\R^{\V(T)}$
with coordinates $(h_i)$. The interior of $\Q(T, \phi)$
corresponds to generalized polyhedra with all edges $\l^*_{ij}$ of
positive length. Function $\vol$ is a homogeneous polynomial of
degree 3 in $(h_i)$, hence differentiable.

\begin{prp} \label{prp:VolDer}
The volume of $Q \in \Q(T, \phi)$ is a smooth function of $h$ with
the following partial derivatives:
\begin{eqnarray}
\frac{\partial \vol}{\partial h_i} & = & F_i \label{eqn:VolHi},\\
\frac{\partial F_i}{\partial h_j} & = & \frac{\l^*_{ij}}{\sin \phi_{ij}} \label{eqn:FiHj},\\
\frac{\partial F_i}{\partial h_i} & = & -\sum_{j \ne i} \cos
\phi_{ij} \frac{\partial F_i}{\partial h_j} \label{eqn:FiHi}.
\end{eqnarray}
\end{prp}
We need the following lemma.
\begin{lem} \label{lem:DerArea}
Face area $F_i$ as a function of variables $h_{ij}$ has the
following partial derivatives:
\begin{equation} \label{eqn:DerArea}
\frac{\partial F_i}{\partial h_{ij}} = \l^*_{ij}.
\end{equation}
\end{lem}
\begin{proof}. It is not hard to see that
$$
\frac{\partial F_i}{\partial h_{ij}} = \frac{1}{2} \left(
\frac{\partial }{\partial h_{ij}}(h_{ij} h_{ijk} + h_{ik} h_{ikj})
+ \frac{\partial }{\partial h_{ij}}(h_{ij} h_{ijl} + h_{ik}
h_{ilj}) \right).
$$
Now the Lemma follows from
\begin{eqnarray*}
\frac{\partial}{\partial h_{ij}} (h_{ij} h_{ijk} + h_{ik} h_{ikj}) & = & h_{ijk} + h_{ij} \frac{\partial h_{ijk}}{\partial h_{ij}} + h_{ik} \frac{\partial h_{ikj}}{\partial h_{ij}}\\
& = & h_{ijk} - h_{ij} \cot \omega_{ijk} + \frac{h_{ik}}{\sin \omega_{ikj}}\\
& = & h_{ijk} + \frac{h_{ik} - h_{ij} \cos \omega{ijk}}{\sin \omega_{ijk}}\\
& = & 2 h_{ijk}\ .
\end{eqnarray*}
\end{proof}

\begin{proof} \textit{of Proposition \ref{prp:VolDer}.}
Equalities (\ref{eqn:FiHj}) and (\ref{eqn:FiHi}) follow from
$$
\frac{\partial F_i}{\partial h_j} = \frac{\partial F_i}{\partial
h_{ij}} \frac{\partial h_{ij}}{\partial h_j}\ ,
$$
respectively from
$$
\frac{\partial F_i}{\partial h_i} = \sum_j \frac{\partial
F_i}{\partial h_{ij}} \frac{\partial h_{ij}}{\partial h_i}\ .
$$
To prove (\ref{eqn:VolHi}), differentiate (\ref{eqn:DVol}) by
Leibniz rule and use formulas (\ref{eqn:FiHj}) and
(\ref{eqn:FiHi}).
\end{proof}

As a byproduct we obtain the following analogue of the Schl\"afli
formula:
\begin{cor}
$$
\sum_i h_i dF_i = 2 \, d \vol\ .
$$
\end{cor}
Note that for classical convex polyhedra formulas
(\ref{eqn:VolHi}) -- (\ref{eqn:DerArea}) have a clear geometric
interpretation.

The triangulation $T$ may have loops and multiple edges. Formulas
(\ref{eqn:FiHj}) and (\ref{eqn:FiHi}) can be generalized to this
case similarly to (\ref{eqn:Complic}).

\subsection{Duality} \label{subsec:Dual}
The formulas of Propositions \ref{prp:CurvDer} and
\ref{prp:VolDer} are closely connected.
\begin{dfn} \label{dfn:DualPol}
Let $P = (T,r)$ be a generalized convex polytope. Its spherical
section centered at the apex produces a generalized simplicial fan
$(T, \phi)$, where $\phi_{ij}$ is the angle of the triangle $aij$
at the vertex $a$. The triple $P^* = (T,\phi,h)$ with $h_i =
\frac{1}{r_i}$ is called the generalized convex polyhedron dual to
$P$.
\end{dfn}
We justify this definition by showing that $\ell^*_{ij} \ge 0$ for
any $ij \in \E(T)$.

\begin{lem}
If $OABC$ is an orthoscheme, then $OC^*B^*A^*$ is also an
orthoscheme. Here $X^*$ is the point on the ray $OX$ such that
$OX^* \cdot OX = 1$.
\end{lem}
\begin{proof}. $OABC$ is an orthoscheme if and only if $OAB, OAC$, and $OBC$ are two-dimensional orthoschemes, that is right triangles. Since triangles $OXY$ and $OY^*X^*$ are similar, $OC^*B^*, OC^*A^*$, and $OB^*A^*$ are also orthoschemes. This implies that $OC^*B^*A^*$ is an orthoscheme. \end{proof}

Call $OC^*B^*A^*$ the orthoscheme dual to $OABC$.
\begin{lem} \label{lem:DualOrth}
Decompose each of the constituting pyramids of $P$ into 6
orthoschemes. Then the duals to these orthoschemes are exactly
those that form the orthoscheme decomposition of $P^*$. The
following formulas hold:
\begin{eqnarray}
h_{ij} & = & \frac{\cot \rho_{ij}}{r_i} \label{eqn:hij*},\\
h_{ijk} & = & \frac{\cot \alpha_{ij}}{r_i \sin\rho_{ij}}
\label{eqn:hijk*}.
\end{eqnarray}
\end{lem}

\begin{figure}[ht]
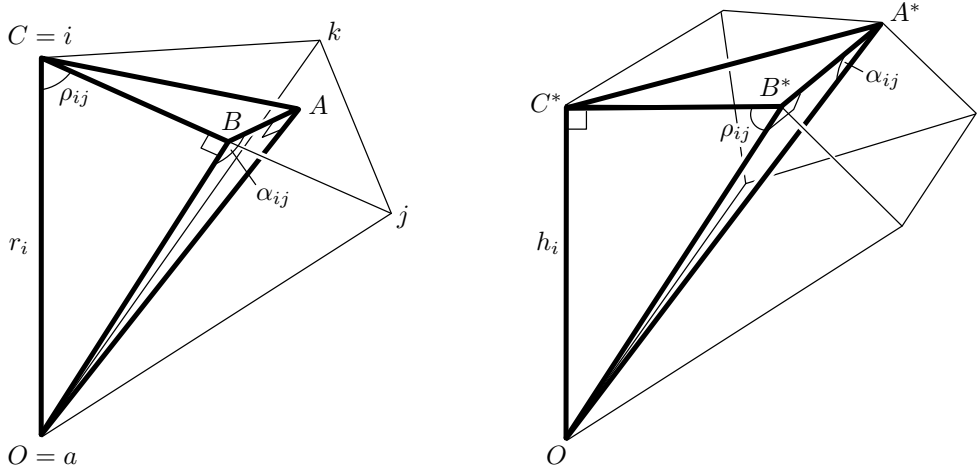

\centerline{\input{Figures/fig06.pstex_t}\hspace{1.7cm}\input{Figures/fig07.pstex_t}}
\caption{Building blocks of a generalized convex polytope and of
its dual polyhedron.} \label{fig:Dual}
\end{figure}

\begin{proof}. For a pyramid $ijk$, the six dual orthoschemes can be identified along common edges and faces. Together they form an object that has combinatorics of the cube, see Figure \ref{fig:Dual}. The three faces of the cube which do not contain the apex are orthogonal to the side edges of the pyramid and have the distances $h_i, h_j, h_k$ from the apex. But exactly this construction was invoked in Definition \ref{dfn:Polyh}. Therefore each of the dual orthoschemes has parameters $h_i, h_{ij}, h_{ijk}$ with indices appropriately permuted.

In the pyramid $ijk$ let $A$ be the foot of the perpendicular to
the base, $B$ be the foot of the perpendicular from $A$ to $ij$,
and $C$ be the vertex $i$. In the dual orthoscheme $C^*B^*A^*$ we
have $A^*B^* = h_{ijk}$ and $B^*C^* = h_{ij}$. To prove the
formulas (\ref{eqn:hij*}) and (\ref{eqn:hijk*}), note that $\angle
OA^*B^* = \angle OBA = \alpha_{ij}$ and $\angle OB^*C^* = \angle
OCB = \rho_{ij}$.
\end{proof}

Substitution of (\ref{eqn:hijk*}) into (\ref{eqn:Dl}) yields
\begin{equation} \label{eqn:l*}
\ell^*_{ij} = \frac{\cot \alpha_{ij} + \cot \alpha_{ji}}{r_i
\sin\rho_{ij}},
\end{equation}
which is non-negative because $\alpha_{ij} + \alpha_{ji} =
\theta_{ij} \le \pi$. This shows that $P^*$ is indeed a
generalized convex polyhedron and justifies Definition
\ref{dfn:DualPol}.

\noindent\emph{Remark.}\quad An orthoscheme can degenerate so that
some of the points $O$, $A$, $B$, $C$ coincide. A formal approach
that includes this possibility is to define an orthoscheme as an
orthonormal basis and a triple of real parameters. Then
degenerations correspond to vanishing of some of the parameters.
If the first parameter vanishes, then so do the others, and the
dual orthoscheme is infinitely large. This does not happen, when
we dualize the orthoschemes in Lemma \ref{lem:DualOrth}, because
the first parameter there is the altitude of a constituting
pyramid and does not vanish ($O \ne A$).

\begin{thm} \label{thm:DerCurv=DerVol}
The Hessian of the total scalar curvature at $P$ is equal to the
Hessian of the volume at $P^*$:
$$
\frac{\partial^2 H}{\partial r_i \partial r_j} (P) =
\frac{\partial^2 \vol}{\partial h_i \partial h_j} (P^*) \quad
\mbox{for all }i,j.
$$
Equivalently, the Jacobian of the curvatures of $P$ with respect
to the radii equals the Jacobian of the face areas of $P^*$ with
respect to the altitudes:
\begin{equation} \label{eqn:Dk=Df}
\frac{\partial \kappa_i}{\partial r_j}(P) = \frac{\partial
F_i}{\partial h_j}(P^*) \quad \mbox{for all }i,j.
\end{equation}
\end{thm}
\begin{proof}. Substitute (\ref{eqn:l*}) into (\ref{eqn:FiHj}) and (\ref{eqn:FiHi}). After applying the sine theorem to triangle $aij$, the right hand sides of (\ref{eqn:FiHj}) and (\ref{eqn:FiHi}) become identical with those of (\ref{eqn:KiRj}) and (\ref{eqn:KiRi}). \end{proof}

\noindent\emph{Remarks.} A generalized simplicial cone $(T,\phi)$
determines a spherical polyhedral metric $S$ on the sphere.
Similarly to Section \ref{sec:CGPandDel}, one can consider the set
$\Q(S)$ of all generalized convex polyhedra with the spherical
section $S$. Then a collection of altitudes $h_i$ defines the
triangulation $T$ uniquely, up to the edges with zero dual edge
length: $\ell^*_{ij} = 0$. Coordinates $(h_i)$ make $\Q(S)$ to a
polyhedron in $\R^n$. The chamber $\Q^T(S)$ that corresponds to a
triangulation $T$ is also a polyhedron and was denoted by
$\Q(T,\phi)$ in Subsection \ref{subsec:GenPolyh}. The polytope $P$
in Theorem \ref{thm:DerCurv=DerVol} has a fixed boundary $M$ and
depends only on $r$; similarly, the polyhedron $P^*$ has a fixed
spherical section $S$ and depends on $h$. Note that the duality
map $P \mapsto P^*$ maps the polytopes with the same boundary to
polyhedra with different spherical sections.

The function $\vol$ is of class $C^2$ on $\Q(S)$, because the
formulas of Proposition \ref{prp:VolDer} agree at the common
boundary points of chambers $\Q^{T_1}(S)$ and $\Q^{T_2}(S)$.

In the classical situation we have the space of polyhedra with
given outer normals $v_i$ to the faces. Minkowski theorem
\cite{Min} says that the condition $\sum F_i v_i = 0$ is necessary
and sufficient for the existence of a polyhedron with the given
face normals and face areas $F_i$. The classical proof of
Minkowski theorem \cite[Section 7.1]{Schn} is based on the
properties of the Hessian of the volume. Because our proof of
Alexandrov's theorem uses non-degeneracy of the Hessian of the
total scalar curvature, Theorem \ref{thm:DerCurv=DerVol} provides
a connection between the Minkowski and Alexandrov's theorems.

\subsection{Mixed volumes of generalized convex polyhedra}
Theorem \ref{thm:DerCurv=DerVol} shows that in order to prove the non-degeneracy of $\left( \frac{\partial \kappa_i}{\partial r_j} \right)$, it suffices to prove the non-degeneracy of  $\left( \frac{\partial^2 \vol}{\partial h_i \partial h_j} \right)$. For classical convex polyhedra the signature of this Hessian is known. This information is expressed as Alexandrov-Fenchel inequalities \cite{Ale37}, and more explicitely in the lemmas used to prove them, see \cite[Section 6.3, Propositions 3 and 4]{Schn}.% Here we show that for generalized polyhedra the Hessian of the volume is non-degenerate; the proof uses methods similar to those in Alexandrov's proof of Alexandrov-Fenchel inequalities, \cite{Ale37}, \cite{Schn}.

Any generalized convex polyhedron has curvatures $\kappa_i$,
defined as the curvatures of its spherical section.
\begin{thm} \label{thm:NonDeg}
Let $Q$ be a generalized convex polyhedron such that
$$
\kappa_i > 0,\, h_i > 0,\, F_i > 0 \quad \mbox{for all }i,
$$
where $\kappa_i$ is the curvature, $h_i$ is the altitude, and
$F_i$ is the area of the $i$-th face. Then the Hessian of the
volume of $Q$ is non-degenerate:
$$
\det \left( \frac{\partial^2 \vol}{\partial h_i \partial h_j}
\right) \ne 0.
$$

\end{thm}
We will prove this theorem using the notion of mixed volumes.

Let us introduce some notation. For any $x \in \R^{\V(T)}$ we
denote by $x_{ij}$ and $x_{ijk}$ linear functions of $x$ obtained
by substituting $x$ for $h$ in formulas (\ref{eqn:Dhij}) and
(\ref{eqn:Dhijk}), respectively. In a similar way we define
functions $\l^*_{ij}(x)$, $F_i(x)$, and $\vol(x)$. Besides,
introduce a function $\overline{F_i}$ on $\R^{\lk i}$ as follows.
The point with coordinates $(x_{ij})_{j \in \lk i}$ in $\R^{\lk
i}$ is sent to $\frac{1}{2} \sum_j x_{ij} \ell^*_{ij}$, where
$\l^*_{ij}$ is again computed via (\ref{eqn:Dl}). In other words,
$$
F_i = \overline{F_i} \circ \pi_i,
$$
where $\pi_i : \R^{\V(T)} \to \R^{\lk i}$ maps $(x_i)$ to
$(x_{ij})$.

It is immediate that all of $F_i$ and $\overline{F_i}$ are
quadratic forms and $\vol$ is a cubic form.
\begin{dfn}
The {\rm mixed areas} $F_i(\cdot,\cdot)$ and
$\overline{F_i}(\cdot,\cdot)$ are the symmetric bilinear forms
associated to the quadratic forms $F_i$ and $\overline{F_i}$. The
{\rm mixed volume} $\vol(\cdot, \cdot, \cdot)$ is the symmetric
trilinear form associated to the cubic form $\vol$.
\end{dfn}

\begin{prp}
The mixed areas and the mixed volume can be computed by the
formulas
\begin{eqnarray}
F_i(x,y) & = & \frac{1}{2} \sum_j x_{ij} \l^*_{ij}(y), \label{eqn:MixArea}\\
\vol(x,y,z) & = & \frac{1}{3} \sum_i x_i F_i(y,z).
\label{eqn:MixVol}
\end{eqnarray}
\end{prp}
\begin{proof}. It suffices to prove (\ref{eqn:MixArea}) for the bilinear form $\overline{F_i}$. For any symmetric bilinear form $B$ we have:
$$
B(x,y) = \frac{1}{2} \sum_i x_i \frac{\partial B(y,y)}{\partial
y_i}.
$$
Thus (\ref{eqn:MixArea}) follows from Lemma \ref{lem:DerArea}.

A similar argument shows that
$$
\vol(x,y,y) = \frac{1}{3} \sum_i x_i F_i(y,y).
$$
Using  $2 \, \vol(x,y,z) = \vol(x,y+z,y+z) - \vol(x,y,y) -
\vol(x,z,z)$ we get (\ref{eqn:MixVol}).
\end{proof}

\begin{lem} \label{lem:FiNonDeg}
If $\kappa_i > 0$, then the form $\overline{F_i}$ is
non-degenerate and has exactly one positive eigenvalue, which is
simple.
\end{lem}
\begin{proof}. First let us prove the non-degeneracy. Assume that for a vector $y \in \R^{\lk i}$ we have $\overline{F_i}(x,y) = 0$ for any $x$. Then from (\ref{eqn:MixArea}) we have $\l^*_{ij}(y) = 0$ for every $j$. We need to show that this implies $y = 0$.

Let us go back to the geometric origins. We have a generalized
polygon with the angles $\omega_{ijk}$ between the normals to
consecutive sides. Numbers $y_{ij}$ are altitudes to the sides. We
can develop the boundary of the polygon onto the plane so that it
becomes a polygonal line. The rotation by the angle $\omega_i =
\sum_{jk} \omega_{ijk}$ sends the starting point of the line to
its endpoint. The equalities $\ell^*_{ij}(y) = 0$ mean that all of
the segments of the line have zero length. Thus the line
degenerates to a point invariant under rotation. Because $\omega_i
= 2\pi - \kappa_i < 2\pi$, this point can be only the origin. Thus
$y_{ij} = 0$ for every $j$.

To determine the signature, deform the collection of angles
$(\omega_{ijk})_{jk}$ continuously so that they become all equal.
If their sum remains less than $2\pi$ during the process of
deformation, then the form $\overline{F_i}$ remains non-degenerate
and it suffices to compute its signature at the final point. This
can be done explicitely by extracting the coefficients from
(\ref{eqn:Dhijk}) and by computing the eigenvalues.

Note that in the classical situation $\omega_i = 2\pi$, and the
polygon can degenerate to any point. This means that the form
$\overline{F_i}$ has  a 2-dimensional nullspace. The form
$\overline{F_i}$ is non-degenerate if and only if $\omega_i$ is
not an integer multiple of $2\pi$. The dimension of the positive
subspace is $2k-1$,for $\omega_i \in (2(k-1)\pi, 2k\pi)$.
\end{proof}

\begin{proof} \emph{of Theorem \ref{thm:NonDeg}.}
We have
$$
\vol(h+x) = \vol(h) + 3\vol(h,h,x) + 3\vol(h,x,x) + \vol(x).
$$
We have to show that under assumption $\kappa_i > 0$ for all $i$,
the quadratic form $\vol(h,\cdot,\cdot)$ is non-degenerate for any
$h$ that satisfies $h_i > 0, F_i(h) > 0$ for all~$i$.

Assume the converse. Then there exists a non-zero vector $x$ such
that $\vol(h,x,y) = 0$ for any $y$. Since $\vol(h,x,y) =
\vol(y,h,x) = \frac{1}{3} \sum_i y_i F_i(h,x)$, this implies
$F_i(h,x) = 0$ for every $i$. Thus $x$ is orthogonal to $h$ with
respect to the form $F_i$. Because $F_i(h,h) > 0$ by assumption,
Lemma \ref{lem:FiNonDeg} implies that $F_i(x,x) \le 0$. Besides,
if $F_i(x,x) = 0$, then $x_{ij} = 0$ for every $j \in \lk i$.

Consider the mixed volume $\vol(h,x,x)$. On one hand we have
$$
\vol(h,x,x) = \frac{1}{3} \sum_i x_i F_i(h,x) = 0.
$$
On the other hand
$$
\vol(h,x,x) = \frac{1}{3} \sum_i h_i F_i(x,x) \le 0,
$$
because $h_i > 0$ for every $i$. Hence $F_i(x,x) = 0$ for every
$i$, and $x_{ij} = 0$ for every edge $ij$. Note that for every $i
\in \V(T)$ we have $\lk i \ne \emptyset$, otherwise $F_i(h) = 0$.
By formula (\ref{eqn:Dhij}), $x_{ij} = 0 = x_{ji}$ together with
$\phi_{ij} \in (0,\pi)$ implies $x_i = 0$. This contradicts the
assumption $x \ne 0$.
\end{proof}

Let $Q = P^*$ be the dual to a generalized convex polytope $P$.
Then $h_i = \frac{1}{r_i} > 0$ for every $i$. Also, $P$ and $P^*$
have identical spherical sections, thus their curvatures are
equal. However, there are examples of positively curved polytopes
$P$ where $F_i(P^*) < 0$ for some $i$. The following Proposition
describes a class of polytopes such that their duals satisfy the
assumptions of Theorem \ref{thm:NonDeg}.

\begin{prp} \label{prp:d>k}
Let $P$ be a generalized convex polytope with boundary $M$ and
suppose that
$$
0 < \kappa_i < \delta_i \quad \mbox{for all }i,
$$
where $\delta_i$ is the angular defect of the $i$-th singularity
on $M$. Then we have
$$
F_i > 0 \quad \mbox{for all }i
$$
in the dual polyhedron $P^*$.
\end{prp}
\begin{proof}. The $i$-th face $P^*_i$ of $P^*$ is a generalized convex polyhedron of dimension 2. It is defined by the angles $(\omega_{ijk})_{jk}$ between the normals to the sides and the altitudes $(h_{ij})_j$ of the sides. Formulas (\ref{eqn:Dhijk}), (\ref{eqn:Dl}), (\ref{eqn:DArea}) allow to compute the side lengths and the area of $P^*_i$. Alternatively, one constructs 2-dimensional orthoschemes, that is right triangles with legs $h_{ij}, h_{ijk}$, and glues them along common sides.

To the Euclidean generalized polyhedron $P^*_i$ there is its
spherical counterpart $SP^*_i$. Take a sphere of radius
$\frac{1}{r_i}$ and place a constituting orthoscheme of $P^*_i$ in
the plane tangent to the north pole so that the base point of the
orthoscheme is at the north pole. Project the orthoscheme to the
northern hemisphere using the sphere center as the center of
projection. The image is a spherical right triangle which we also
call an orthoscheme, see Figure \ref{fig:OrthProj}. The polyhedron
$SP^*_i$ is defined as the complex of thus constructed spherical
orthoschemes. Because of (\ref{eqn:hij*}), the first parameter of
the spherical orthoscheme equals $\frac{\pi}{2} - \rho_{ij}$.
Therefore $SP^*_i$ is defined by the angles $\omega_{ijk}$ between
the normals to the sides and the altitudes $\frac{\pi}{2} -
\rho_{ij}$.

\begin{figure}[ht]
\centerline{\input{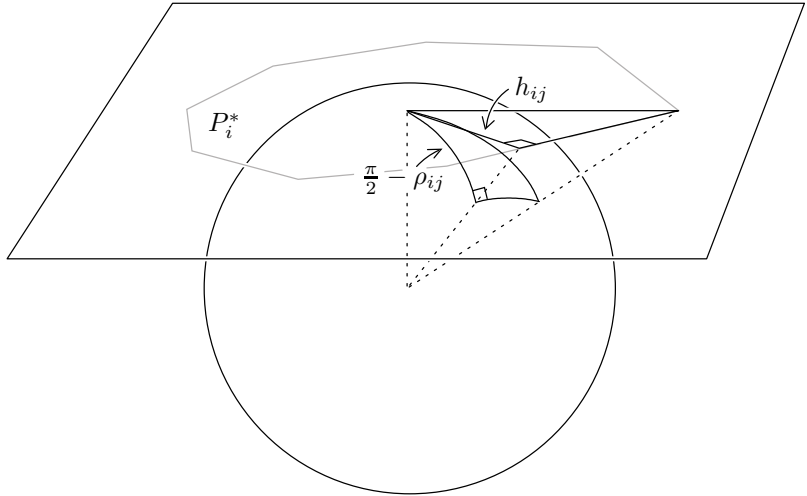}} \caption{To the
definition of $SP^*_i$.} \label{fig:OrthProj}
\end{figure}

We claim that
\begin{equation} \label{eqn:AreaDK}
\area(SP^*_i) = \delta_i - \kappa_i.
\end{equation}
By Lemma \ref{lem:Sph5} this together with the assumption
$\delta_i > \kappa_i$ implies $F_i > 0$.

Let $SP_i$ be the spherical section of $P$ at the vertex $i$. It
consists of the spherical triangles with two sides of length
$\rho_{ij}, \rho_{ik}$ and angle $\omega_{ijk}$ between them.
Denote by $\gamma_{ijk}$ the length of the third side of this
triangle. The spherical generalized 2-polytop $SP_i$ has $SP^*_i$
as its spherical dual, see Figure \ref{fig:SP}. Consider the pair
of orthoschemes $(ij, ijk)$ and $(ik, ikj)$ in $SP^*_i$. Together
they form a quadrilateral with two right opposite angles. The
angle of this quadrilateral at the center of the polyhedron equals
$\omega_{ijk}$, whereas the fourth angle is $\pi - \gamma_{ijk}$.
Therefore the area of the quadrilateral is $\omega_{ijk} -
\gamma_{ijk}$. This remains true when some parameters of the
orthoschemes are negative, if one considers signed angles and
signed areas. As a result, we have
$$
\area(SP^*_i) = \sum_{jk} (\omega_{ijk} - \gamma_{ijk}) = \omega_i
- \per(SP_i).
$$
Since $\kappa_i = 2\pi - \omega_i$ and $\delta_i = 2\pi -
\per(SP_i)$, this implies the equality (\ref{eqn:AreaDK}).
\end{proof}

\begin{figure}[ht]
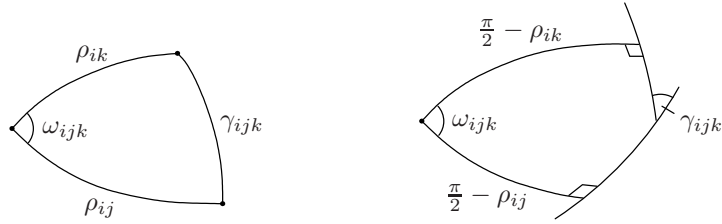

\centerline{\input{Figures/fig09.pstex_t}\hspace{2.5cm}\input{Figures/fig10.pstex_t}}
\caption{The spherical generalized polygons $SP_i$ and $SP^*_i$.}
\label{fig:SP}
\end{figure}

For a 2-dimensional generalized convex polyhedron $Q$ define its
\emph{spherical image} $SQ$ exactly as we defined $SP_i^*$ in the
proof of Proposition \ref{prp:d>k}, but taking a unit sphere
instead, see Figure \ref{fig:OrthProj}.
\begin{lem} \label{lem:Sph5}
Let $C$ be a positively curved generalized convex polyhedron of
dimension 2, and let $SC$ be its spherical image. Suppose that
$\area(SC) > 0$. Then also $\area(C) > 0$.
\end{lem}
\begin{proof}. If all of the altitudes of $C$ are negative, then so are the altitudes of $SC$, which leads to $\area(SC) < 0$.

Let $h_i$ be a positive altitude of $C$. We have $h_{ij} + h_{ik}
= \l^*_i \ge 0$, where $\l^*_i$ is the length of the $i$-th side
of $C$, and $h_{ij}, h_{ik}$ are the signed lengths of the
segments into which the side is divided by the foot of the
altitude. Without loss of generality $h_{ij} \ge 0$. Cut $C$ along
the segment joining the vertex $ij$ to the apex. Develop the
result onto the plane. The two sides of the cut become sides of an
isosceles triangle $\Delta$, and the boundary of $C$ developes to
a convex polygonal line joining the two base vertices of $\Delta$.
It can be shown that, when extended by the base of $\Delta$, the
line remains convex. Let $D$ be the convex polygon bounded by it.
If $\kappa \ge \pi$, then $\Delta$ and $D$ lie on different sides
from their common edge, and one easily sees that all of the
altitudes of $C$ are positive. This automatically implies
$\area(C)>0$. If $\kappa < \pi$, then we have
$$
\area(C) = \area(D) - \area(\Delta),
$$
and, if $J$ denotes the Jacobian of the projection to the sphere,
$$
\area(SC) = \int_D J d\area - \int_\Delta J d\area.
$$

\begin{figure}[ht]
\centerline{\input{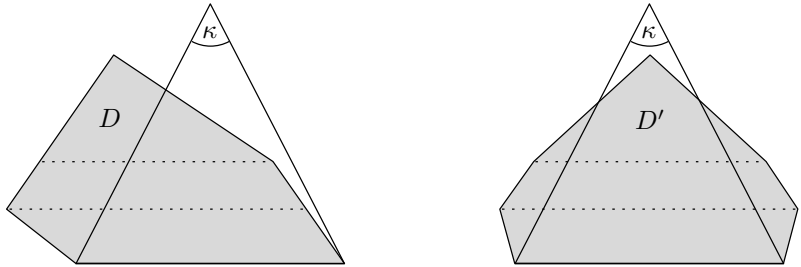}} \caption{Triangle
$\Delta$, polygon $D$ and its Steiner symmetral $D'$.}
\end{figure}

To show that $\area(SC) > 0$ implies $\area(C) > 0$, replace $D$
by its Steiner symmetral $D'$ with respect to the line
perpendicular to the base of $\Delta$. This means that $D$ gets
decomposed into segments parallel to the base and every segment
gets translated so that its midpoint lies on the altitude of
$\Delta$. The polygon $D'$ is symmetric with respect to the
altitude of $\Delta$ and has the same area as $D$. On the other
hand, it is not hard to see that $\int_{D'} J d\area \ge \int_D J
d\area$. This implies
$$
\int_{D'} J d\area > \int_\Delta J d\area.
$$
It is easy to see that for any $x \in D' \setminus \Delta$ and for
any $y \in \Delta \setminus D'$ we have $J(x) \le J(y)$. It
follows that $\area(D') > \area(\Delta)$ and thus $\area(C) > 0$.
\end{proof}

\begin{cor} \label{cor:NonDeg}
Let $P$ be a generalized convex polytope with radii $(r_i)$ and
curvatures $(\kappa_i)$. Suppose that
$$
0 < \kappa_i < \delta_i \quad \mbox{for all }i,
$$
where $\delta_i$ is the angular defect of the $i$th singularity on
the boundary of $P$. Then the matrix $\left( \frac{\partial
\kappa_i}{\partial r_j} \right)$ is non-degenerate.
\end{cor}
\begin{proof}. This follows from Proposition \ref{prp:d>k}, Theorem \ref{thm:NonDeg}, and from the fact that the Jacobian of the map $r \mapsto \kappa$ equals the Hessian of $H$, see (\ref{eqn:DmDr}). \end{proof}

Now we study the matrix $\left( \frac{\partial \kappa_i}{\partial
r_j} \right)$ at $\kappa_i = 0$ for all $i$, that is in the case
when the generalized convex polytope is a usual convex polytope.
As a corollary we get a proof of the infinitesimal rigidity of
polytopes.

\begin{prp} \label{prp:KerDescr}
Let $P$ be a generalized convex polytope with all curvatures zero.
Then the Jacobian $\left( \frac{\partial \kappa_i}{\partial r_j}
\right)$ has corank $3$. A vector $dr = (dr_i)$ lies in the kernel
of $\left( \frac{\partial \kappa_i}{\partial r_j} \right)$ if and
only if there exists $x \in \R^3$ such that
\begin{equation} \label{eqn:ker}
dr_i = \left\langle \frac{p_i-a}{\|p_i - a\|}, x \right\rangle
\quad \mbox{for all }i.
\end{equation}
Here $p_i$ are the vertices of $P$, and $a$ is the apex, for some
realization of $P$ as a convex polytope in $\R^3$.
\end{prp}
\begin{proof}. Use the equality (\ref{eqn:Dk=Df}). Here $P^*$ is the convex polytope in $\R^3$ polar to $P$. It is known that the matrix $\left( \frac{\partial F_i}{\partial h_j} \right)$ has corank $3$, see \cite[Section 6.3, Proposition 3]{Schn}. Hence the Jacobian $\left( \frac{\partial \kappa_i}{\partial r_j} \right)$ has also corank $3$.

Now let us find vectors in the kernel. Let $x \in \R^3$ be any
vector. Consider the family $\{P_x(t)|\: t \in (-\epsilon,
\epsilon)\}$ of generalized polytopes with zero curvatures, where
$P_x(t) \subset \R^3$ has vertices $p_i$ and the apex $a-tx$. From
$r_i = \|p_i - a\|$ it is easy to get
$$
\left.\frac{dr_i}{dt}\right|_{t=0} = \left\langle
\frac{p_i-a}{\|p_i - a\|}, x \right\rangle.
$$
Therefore any infinitesimal deformation given by (\ref{eqn:ker})
results in $d\kappa_i = 0$ for all $i$.
\end{proof}

\begin{cor}
Convex polytopes are infinitesimally rigid.
\end{cor}
\begin{proof}. Let $P \subset \R^3$ be a convex polytope with vertices $(p_i)$. Triangulate its non-triangular faces by diagonals. Let $(dp_i)$ be infinitesimal motions of the vertices that result in a zero first order variation of the edge lengths. Take any point $a$ in the interior of $P$ and cut $P$ into pyramids with the apex $a$. The deformation changes the dihedral angles of the pyramids, but the total angles $\omega_i$ around the radial edges remain equal to $2\pi$. This means that the vector $(dr_i = d\|p_i - a\|)$ belongs to the kernel of the Jacobian $\left( \frac{\partial \omega_i}{\partial r_j} \right)$. Recall that $\kappa_i = 2\pi - \omega_i$. By Proposition \ref{prp:KerDescr} we have (\ref{eqn:ker}). It follows that $(dp_i)$ is the restriction of the differential of a rigid motion of $P$.
\end{proof}

The infinitesimal rigidity of convex polytopes is a classical
result due to Dehn \cite{Dehn}. Recently different new proofs were
suggested by Filliman \cite{Filliman}, Schlenker \cite{Sch1} and
Pak \cite{Pak}.

The proof of Filliman \cite{Filliman} is very similar to ours and
is formulated in the language of stresses on frameworks.

Schlenker in \cite{Sch1} studies infinitesimal isometric
deformations of spherical and Euclidean polygons and the
corresponding first-order variations of their angles. He
introduces a quadratic form on the space of infinitesimal
deformations and proves a remarkable positivity property for it.

\section{Proof of Alexandrov's theorem} \label{sec:Proof}
In this section we prove the existence part of Theorem
\ref{thm:Alex}. The uniqueness can be proved by extending the
classical Cauchy's argument, see \cite[Section 3.3]{Al05}. In
order to establish existence, we construct a family $\{P(t)\}_{t
\in [0,1]}$ of generalized convex polytopes with the boundary $M$
such that the curvatures of $P(0)$ vanish. Then $P(0)$ is a
desired convex polytope.

\begin{thm} \label{thm:Family}
For any convex Euclidean polyhedral metric $M$ on the sphere there
exists a family $\{P(t) = (T(t),r(t))|\: t \in (0,1]\}$ of
generalized convex polytopes with boundary $M$ such that
\begin{enumerate}
\item \label{ittt:1} $T(1)$ is the Delaunay triangulation of $M$,
and $r_i(1) = R$ for all $i$ and a sufficiently large $R$; \item
\label{ittt:2} the curvatures $\kappa_i(t)$ of the polytope $P(t)$
are proportional to $t$:
\begin{equation} \label{eqn:kappa(t)}
\kappa_i(t) = t \cdot \kappa_i(1) \quad \mbox{ for all } i.
\end{equation}
\item \label{ittt:3} $r(t)$ is of class $C^1$ on $(0,1]$; \item
\label{ittt:4} there exists the limit $r(0) = \lim_{t \to 0}
r(t)$; \item \label{ittt:5} there exists a convex polytope $P
\subset \R^3$ with boundary isometric to $M$ and vertices $p_i$
such that $r_i(0) = \|p_i - a\|$, where $a \in P$ is the unique
point that satisfies the condition
\begin{equation} \label{eqn:apex}
\sum_i \kappa_i(1) \frac{p_i - a}{\|p_i - a\|} = 0.
\end{equation}
\end{enumerate}
\end{thm}
Theorem \ref{thm:Family} provides the following recipe for
construction of a polytope with the boundary $M$. The path
$P(t)_{t \in (0,1]}$ in the space $\P(M)$ is the preimage of the
path $\kappa(t)_{t \in [0,1]}$ given by (\ref{eqn:kappa(t)}).
Since a generalized convex polytope is uniquely determined by its
radii, we have a system of differential equations:
\begin{equation} \label{eqn:System}
\frac{dr}{dt} = \left(\frac{\partial \kappa_i}{\partial
r_j}\right)^{-1} \kappa(1).
\end{equation}
Theorem \ref{thm:Family} states that the system has a solution on
$(0,1]$ with the initial condition $r(1) = (R,\ldots,R)$ for $R$
sufficiently large. Besides, $r(t)$ converges at $t=0$. At this
point there are two possibilities: either the pyramids with edge
lengths $r_i(0)$ over some triangulation $T(0)$ form a convex
polytope, or all of the pyramids degenerate and there bases form a
doubly-covered polygon. When solving the system (\ref{eqn:System})
numerically one should not forget that the formulas
(\ref{eqn:KiRj}), (\ref{eqn:KiRi}) for the partial derivatives
$\frac{\partial \kappa_i}{\partial r_j}$ depend on the
combinatorics of the triangulation $T(t)$. During the deformation
of the generalized convex polytope $P(t)$ the triangulation $T(t)$
may change. In a generic case there are flips that happen at
different moments.

\subsection{Proof of Theorem \ref{thm:Family}}
The proof uses the description of the space $\P(M)$ of generalized
convex polytopes and the non-degeneracy of the Hessian of
functional $H$. Especially we refer to Proposition \ref{prp:P(M)}
and Corollary \ref{cor:NonDeg}. Several lemmas needed to deal with
degenerating generalized polytopes are postponed to Section
\ref{sec:Lemmas}.

In Section \ref{sec:CGPandDel} we used an embedding of $\P(M)$
into $\R^n$ via squares of the radii. Here we consider radii as
coordinates on $\P(M)$. Thus a generalized convex polytope is
identified with a point $r = (r_1,\ldots,r_n) \in \R^n$. The map
\begin{eqnarray}
\P(M) & \to & \R^n, \label{eqn:map}\\
r & \mapsto & \kappa \nonumber,
\end{eqnarray}
is of class $C^1$ in the interior of $\P(M)$ by Proposition
\ref{prp:CurvDer}. We have to show that the path $\kappa:(0,1] \to
\R^n$ given by (\ref{eqn:kappa(t)}) can be lifted through the map
(\ref{eqn:map}) to a path $r:(0,1] \to \P(M)$. The question of
convergence of $r$ at $t \to 0$ will be treated later.

\begin{lem} \label{lem:Start}
For $R$ large enough, there is a lift of the path
(\ref{eqn:kappa(t)}) on an interval $(t_0,1]$ for some $t_0$.
\end{lem}
\begin{proof}. By property \ref{itt:1}. from Proposition \ref{prp:P(M)}, there is a generalized convex polytope $P(1)$ with large equal radii. Moreover, $P(1)$ lies in the interior of the space $\P(M)$. The corresponding triangulation $T(1)$ of $M$ is the weighted Delaunay triangulation with equal weights, that is the Delaunay triangulation.

The spherical section of $P(1)$ at the vertex $i$ is a convex star
polygon, see Figure \ref{fig:fig12}. We have $\rho_j <
\frac{\pi}{2}$ for every $j$, because $\rho_j$ is an angle at the
base of an isosceles triangle. Therefore by Lemma \ref{lem:Sph1}
we have $\kappa_i(1) < \delta_i$. Besides, as $R$ tends to
infinity, every $\rho_j$ tends to $\frac{\pi}{2}$. This implies
that $\kappa_i(1)$ is close to $\delta_i$ for $R$ sufficiently
large; in particular $\kappa_i(1) > 0$. Thus we have
\begin{equation} \label{eqn:0kd}
0 < \kappa_i(1) < \delta_i \quad \mbox{ for all }i.
\end{equation}
By Corollary \ref{cor:NonDeg}, this implies that the Jacobian of
(\ref{eqn:map}) at $r(1)$ is non-degenerate. The Lemma follows
from the inverse function theorem.
\end{proof}

\begin{lem} \label{lem:open}
If $r:(t_0,1] \to \P(M)$ and $r':(t'_0,1] \to \P(M)$ are two lifts
of the path (\ref{eqn:kappa(t)}), then they coincide on $(t_0,1]
\cap (t'_0,1]$.

The maximum interval where the lift exists is open in $(0,1]$.
\end{lem}
\begin{proof}. Assume the converse to the first statement of the lemma. Let
$$
t_1 = \inf\{t|\: r(t) = r'(t)\}.
$$
Then we have $r(t_1) = r'(t_1)$. By (\ref{eqn:kappa(t)}) and
(\ref{eqn:0kd}), we have
$$
0 < \kappa_i(t_1) < \delta_i \quad \mbox{ for all }i.
$$
Let us show that $r(t_1)$ is an interior point of $\P(M)$. If
$r(t_1)$ is a boundary point of $\P(M)$, then by Proposition
\ref{prp:P(M)}, $r(t_1)$ is also a boundary point of $\D(M)$. This
means that one of the inequalities (\ref{eqn:ConcQuadr}) or
(\ref{eqn:CloGeod}) becomes an equality after the substitution of
$r^2(t_1)$ for $q$. If we have equality in (\ref{eqn:ConcQuadr}),
then the spherical section at the vertex $j$ contains a boundary
geodesic arc of length at least $\pi$. This contradicts Lemma
\ref{lem:Sph2}. If an equality in (\ref{eqn:CloGeod}) is achieved,
then in the pyramid over the triangle $iji$ the vertex $j$ is the
foot of the altitude. This implies $\kappa_j(t_1) = \delta_j$ that
contradicts to (\ref{eqn:0kd}).

By Corollary \ref{cor:NonDeg} and inverse function theorem, the
lift is unique in a neighborhood of $t_1$. This implies the first
statement of the lemma.

The same argument shows that a lift on any closed interval
$[t_1,1] \subset (0,1]$ can be extended to a neighborhood of
$t_1$. This implies the second statement of the lemma. \end{proof}

Let $(t_0,1] \subset (0,1]$ be the maximum interval where the lift
of (\ref{eqn:kappa(t)}) exists. We have to prove $t_0=0$.

\begin{lem} \label{lem:seq}
Let $(t_0,1] \subset (0,1]$ be the maximum interval where the lift
$r: (t_0,1] \to \P(M)$ exists, and let $t_0 > 0$. Assume that
there is a sequence $(t_n)$ in $(t_0,1]$ such that $t_n$ converges
to $t_0$ and $r(t_n)$ converges to a point $r$ in $\P(M)$. Then
the path $r: (t_0,1] \to \P(M)$ can be extended to $t_0$ in a
$C^1$ way by putting $r(t_0) = r$.
\end{lem}
\begin{proof}. Let $\kappa_i$ be the curvatures at the point $r = \lim_{n \to \infty} r(t_n)$. Since $\kappa_i(t_n) = t_n \cdot \kappa_i(1)$ and due to the continuity of the map (\ref{eqn:map}), we have $\kappa_i = t_0 \cdot \kappa_i(1)$. Therefore $0 < \kappa_i < \delta_i$. As in the proof of Lemma \ref{lem:open}, this implies that $r$ is an interior point of $\P(M)$. Besides, by Corollary \ref{cor:NonDeg}, the Jacobian of the map (\ref{eqn:map}) at $r$ is non-degenerate. Thus (\ref{eqn:map}) maps a neigborhood $U$ of $r$ diffeomorphically to a neighborhood $V$ of $\kappa$. For $n$ large enough, $r(t_n)$ lies in $U$. Since the path $\kappa(t)$ is regular at $t = t_0$, the path $r(t)$ extends smoothly to a neighborhood of $t_0$. \end{proof}

In order to show $t_0 = 0$ it suffices to find a sequence $(t_n)$ as in Lemma \ref{lem:seq}. It is not hard to satisfy the condition that $r(t_n)$ converges: one needs only boundedness of $r$ on $(t_0,1]$ which is proved in Lemma \ref{lem:Bound} below. A major problem is to show that the limit lies in $\P(M)$, Lemmas \ref{lem:xFace}--\ref{lem:xVert}. % The converse would mean that the sequence $P(t_n)$ tends to a degenerate generalized polytope. It requires some work to exclude all possible degenerations, and we do it in Lemmas -- . Along the way we use several lemmas from Subsection \ref{subsec:Lemmas}.
\begin{dfn} \label{dfn:SphSec}
The spherical section of a generalized convex polytope $P$ is a
sphere with a spherical polyhedral metric that is obtained by
gluing together the sections of constituting pyramids of $P$ with
the center at the apex.
\end{dfn}
Clearly, the curvatures of the spherical section equal the
curvatures of the polytope. The area of the spherical section can
be computed as
\begin{equation} \label{eqn:SphSecArea}
\area(S) = 4\pi - \sum_i \kappa_i.
\end{equation}
For a generalized convex polytope $P(t)$ denote its spherical
section by $S(t)$. Instead of $P(t_n)$ and $S(t_n)$ we write
briefly $P_n$ and $S_n$.

\begin{lem} \label{lem:Bound}
Let $r: (t_0,1] \to \P(M)$ be a lift of the path
(\ref{eqn:kappa(t)}) through the map (\ref{eqn:map}) for some $t_0
\ge 0$. Then $r$ is bounded on $(t_0,1]$.
\end{lem}
\begin{proof}. Assume that $r_i(t) > C$ for some $t \in (t_0,1]$, some $i$, and some large $C$. Then the triangle inequality in the metric space $P(t)$ implies that $\dist(0,x) > C - \diam(M)$ for any $x \in M$, where $0$ stands for the apex of $P(t)$. We claim that
\begin{equation} \label{eqn:Area<}
\area(S(t)) < \frac{\area(M)}{(C - \diam(M))^2}.
\end{equation}
For $C$ large enough this contradicts to
$$
\area(S(t)) = 4\pi - \sum_i \kappa_i(t) > 4\pi - \sum_i
\kappa_i(1) = \area(S(1)).
$$
In order to prove (\ref{eqn:Area<}), consider a constituting
pyramid of $P(t)$. Denote by $\alpha$ its solid angle at the apex,
and by $A$ the area of its base. By our assumption, the distance
from the apex to any point of the base is larger than $C -
\diam(M)$. It follows that the Jacobian of the central projection
of the base to the unit sphere centered at the apex is less than
$(C - \diam(M))^{-2}$. Thus we have $\alpha < A(C -
\diam(M))^{-2}$. This implies the inequality (\ref{eqn:Area<}),
and the lemma is proved.
\end{proof}

It follows that there is a sequence $(t_n)$ such that $r(t_n)$
converges to some point $r \in \R^n$. Let $\widetilde{q_n}$ be the
$Q$-concave $PQ$-function on $M$ that takes value $r_i^2$ at the
singularity $i$, for all $i$. Recall that $\widetilde{q_n}$ is the
squared distance from the apex of the polytope $P_n$ to its
boundary $M$, and that the space $\P(M)$ is identified with the
space of positive $Q$-concave $PQ$-functions on $M$, see Lemma
\ref{lem:qTr_prop}. By property \ref{it:5}. from Proposition
\ref{prp:PropDel}, $\widetilde{q_n}$ depends continuously on
$r(t_n)$ in the $L^{\infty}$-metric. Since the space $\D(M)$ of
all $Q$-concave $PQ$-functions is closed, this implies that the
sequence $(\widetilde{q_n})$ converges uniformly to a $Q$-concave
$PQ$-function $\widetilde{q}$. Function $\widetilde{q}$ takes
value $r_i$ at the singularity $i$, for all $i$, and is
non-negative. If $\widetilde{q}$ is positive, then the point $r$
lies in $\P(M)$, and we are done. Otherwise there is a point $x
\in M$ such that $\widetilde{q}(x)=0$. Geometrically this means
that in the polytopes $P_n$ the apex tends to the point $x \in M =
\partial P_n$. Since the decomposition $\P(M) = \cup_T \P^T(M)$ is
locally finite by Proposition \ref{prp:P(M)}, we can assume,
replacing $(t_n)$ by its subsequence when needed, that all of the
polytopes $P_n$ have the same triangulation $T(t_n) = T$. Now
there are three possibilities: $x$ lies in the interior of a face
of $T$, or $x$ lies in the interior of an edge of $T$, or $x$ is a
singularity $i$ of $M$.

\begin{lem} \label{lem:xFace}
The point $x$ cannot be an interior point of a face of $T$.
\end{lem}
\begin{proof}. Assume the converse. Let $\Delta \in \F(T)$ be the triangle that contains $x$ in the interior, and denote by $S\Delta_n$ the corresponding spherical triangle in $S_n$. As $n \to \infty$, the spherical triangle $S\Delta_n$ tends to a hemisphere, and its angles tend to $\pi$. There are no self-identifications on the boundary of $S\Delta_n$, because $\kappa_i(t_n) \to \kappa_i(t_0) > 0$. Therefore the complement $S\Delta'_n$ to $S\Delta_n$ in $S_n$ is a singular spherical triangle in the sense of Definition \ref{dfn:SingSphPol}. If $i$ is a vertex of $\Delta$, then the angle of $S\Delta'_n$ at $i$ tends to $\pi - \kappa_i(t_0)$. The side lengths of $S\Delta'_n$ also have limits as $n\to \infty$, these are the angles under which the sides of the triangle $\Delta$ are seen from the point $x$. On the other hand, the perimeter of $S\Delta'_n$ equals the perimeter of $S\Delta_n$ that tends to $2\pi$. Therefore, for sufficiently large $n$ we get a contradiction to Lemma \ref{lem:Sph3}. \end{proof}

\begin{lem} \label{lem:xEdge}
The point $x$ cannot be an interior point of an edge of $T$.
\end{lem}
\begin{proof}. Assume that $x$ is an interior point of an edge $ij \in \E(T)$. Let $\Delta^1, \Delta^2 \in \F(T)$ be the triangles incident to $ij$, and let $S\Delta^1_n$, $S\Delta^2_n$ be their spherical images in the spherical section of the polytope $P_n$. Consider the spherical quadrilateral $S\Diamond_n = S\Delta^1_n \cap S\Delta^2_n$. As $n$ goes to $\infty$, the length of the diagonal $ij$ in $S\Diamond_n$ tends to $\pi$. At the same time, the side lengths of $S\Diamond_n$ have limits different from $0$ and $\pi$. The angles of $S\Diamond_n$ opposite to the diagonal $ij$ tend to $\pi$, whereas the angles $\alpha_n$ and $\beta_n$ at the vertices $i$ and $j$ may have no limits. However, the convexity of $P_n$ implies that
$$
\liminf \alpha_n \ge \pi \quad \mbox{and} \quad \liminf \beta_n
\ge \pi.
$$
It follows that there are no self-identifications on the boundary
of $S\Diamond_n$, thus its complement $S\Diamond'_n$ is a singular
spherical quadrilateral. The perimeter of $S\Diamond'_n$ tends to
$2\pi$, whereas its side lengths have limits different from $0$
and $\pi$. Also, due to
$$
\lim_{n \to \infty} \kappa_k(t_n) = t_0 \cdot \kappa_k(1) > 0
\quad \mbox{for all } k,$$ the angles of $S\Diamond'_n$ have upper
limits less than $\pi$. This contradicts Lemma \ref{lem:Sph3}.
\end{proof}

The most delicate case is when $x$ is a singularity of $M$. To
deal with it we need the hard Lemma \ref{lem:Sph4}.

\begin{figure}[ht]
\centerline{\input{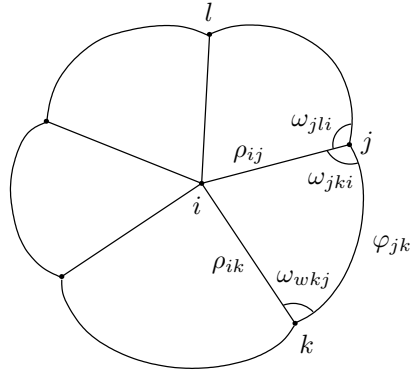}} \caption{The spherical
image of the star of vertex $i$ as $r_i$ tends to $0$.}
\label{fig:fig19}
\end{figure}

\begin{lem}
Function $\widetilde{q}$ cannot vanish at more than one vertex of
$T$.
\end{lem}
\begin{proof}. Assume $\widetilde{q}(i) = 0$, that is $r_i = 0$. Consider the star $\mathrm{st} i$ of the vertex $i$ in the triangulation $T$ and its spherical image $S (\mathrm{st} i)_n$ in $S_n$, see Figure \ref{fig:fig19}. Here all of the angles and lengths depend on $n$. Denote by
$$
\eta_j(n) = \omega_{jki}(n) + \omega_{jli}(n)
$$
the angle at the vertex $j$ in $S (\mathrm{st} i)_n$. We have
$$
\area(S (\mathrm{st} i)_n) = 2\pi - \kappa_i(t_n) - \sum_{j \in
\lk i} (\pi - \eta_j(n)).
$$
By convexity of $P_n$,
\begin{equation} \label{eqn:eta}
\liminf_{n \to \infty} \eta_j(n) \ge \pi.
\end{equation}
Therefore we have
$$
\liminf_{n \to \infty} \area(S (\mathrm{st} i)_n) \ge 2\pi - t_0
\cdot \kappa_i(1).
$$
On the other hand, by (\ref{eqn:SphSecArea})
$$
\lim_{n \to \infty} \area(S_n) = 4\pi - t_0 \cdot \sum_i
\kappa_i(1).
$$
Thus by the first part of Lemma \ref{lem:Sph4} we have
\begin{equation} \label{eqn:Half}
\liminf_{n \to \infty} \area(S (\mathrm{st} i)_n) > \frac{1}{2}
\lim_{n \to \infty} \area(S_n).
\end{equation}
Now, if there is another singularity $m$ such that $r_m \to 0$,
then the stars of $i$ and $m$ have disjoint interiors. On the
other hand, by (\ref{eqn:Half}) both $S (\mathrm{st} i)_n$ and $S
(\mathrm{st} m)_n$ make at least a half of $S_n$, for $n$ large
enough. This contradiction proves the lemma. \end{proof}

\begin{lem} \label{lem:xVert}
Function $\widetilde{q}$ cannot vanish at a vertex of $T$.
\end{lem}
\begin{proof}. It is easy to see that
\begin{equation} \label{eqn:GammaLim}
\lim_{n \to \infty} \phi_{jk}(n) = \gamma_{ijk},
\end{equation}
where $\gamma_{ijk}$ is the angle at the vertex $i$ in the
Euclidean triangle $ijk \in \F(T)$. Consider the angles
$\omega_{kji}(n), \omega_{jki}(n), \omega_{jli}(n)$ and so on, and
call them the {\em base angles} of the star. The base angles can
behave badly. By choosing a subsequence if needed, we can assume
that there exist limits
$$
\lim_{n \to \infty} \omega_{jki}(n) = \omega_{jki} \in [0,\pi]
\quad \mbox{etc.}
$$

{\em Case 1.} The star of $i$ does not degenerate, that is the
limits of all of the base angles are in $(0,\pi)$. In this case
the sequence of spherical sections $(S_n)$ converge to a spherical
polyhedral metric that satisfies all of the assumptions of Lemma
\ref{lem:Sph4}. Here the distinguished singularity $0$ is the
vertex $i$, and the star polygon $C$ is the limit of the stars $S
(\mathrm{st} i)_n$. Angles of $C$ are greater or equal $\pi$
because of the convexity of the polytope $P_n$. Note that the
curvatures of the resulting polyhedral metric equal $t_0 \cdot
\kappa_i(1)$ for all $i \in \V(T)$. Thus by (\ref{eqn:Last}) we
have
\begin{equation} \label{eqn:ohoho}
\kappa_i(1) \ge \left( 1 - \frac{\per(C)}{2\pi} \right) \, \sum_{j
\ne i} \kappa_j(1).
\end{equation}
Because of (\ref{eqn:GammaLim}) and (\ref{eqn:0kd}) we have
$$
\per(C) = 2\pi - \delta_i < 2\pi - \kappa_i(1).
$$
Therefore (\ref{eqn:ohoho}) implies
\begin{equation} \label{eqn:uhuhu}
\sum_{j \ne i} \kappa_j(1) < 2\pi.
\end{equation}
On the other hand, we have
$$
\sum_{j \ne i} \delta_j > 2\pi,
$$
since $\sum_{j \in \Sigma} \delta_j = 4\pi$ and $\delta_i < 2\pi$.
But, if the number $R$ was chosen sufficiently large, then
$\kappa_j(1)$ are close to $\delta_j$ (see the proof of Lemma
\ref{lem:Start}) and also satisfy $\sum_{j \ne i} \kappa_j(1) >
2\pi$. This contradicts the equality (\ref{eqn:uhuhu}), and Case
1. is excluded.

Let us show that if some of the base angles tends to $0$ or $\pi$,
then we can assume without loss of generality
\begin{equation} \label{eqn:wlog}
\omega_{kji} = \pi, \quad \omega_{jki} > 0.
\end{equation}
Let one of the base angles tend to $\pi$, say $\omega_{kji} =
\pi$. If (\ref{eqn:wlog}) fails, then $\omega_{jki} = 0$. Then
(\ref{eqn:eta}) implies $\omega_{jli} = \pi$. If $\omega_{lji} >
0$, then after relabeling we have (\ref{eqn:wlog}). So again
$\omega_{lji} = 0$. Proceeding in this manner we show that all of
the base angles tend to $0$ or $\pi$ depending on their
orientation. But this implies
$$
\rho_{ik}(n) < \rho_{ij}(n) < \rho_{il}(n) < \cdots <
\rho_{ik}(n),
$$
for a sufficiently large $n$. We obtain a similar chain of
inequalities if we assume that one of the base angles tends to
$0$, and $\ref{eqn:wlog}$ fails for any edge $ij \in \lk i$. This
contradiction shows that we may assume (\ref{eqn:wlog}).

{\em Case 2.} We have $\omega_{kji} = \pi, \omega_{jki} \in
(0,\pi)$. This implies
$$
\lim_{n \to \infty} \rho_{ij}(n) = \pi, \quad \lim_{n \to \infty}
\omega_{lji} = \pi.
$$
To the sqherical quadrilateral $ikjl$ the same argument as in the
proof of Lemma \ref{lem:xEdge} can be applied and leads to a
contradiction.

{\em Case 3.} We have $\omega_{kji} = \omega_{jki} = \pi$. We may
assume that the limits
$$
\lim_{n \to \infty} \rho_{ij}(n) = \rho_{ij}, \quad \lim_{n \to
\infty} \rho_{ik}(n) = \rho_{ik}
$$
exist. If both are less than $\pi$, then we apply the argument
from Lemma \ref{lem:xFace} to the spherical triangle $ijk$. If
$\rho_{ij} = \pi$, then both $\rho_{ik}$ and $\rho_{il}$ are less
than $\pi$. Then we consider the spherical quadrilateral $ikjl$
and apply to it the argument of Lemma \ref{lem:xEdge}.
\end{proof}

We proved that $t_0 = 0$ for the maximum lift $r:(t_0,1] \to
\P(M)$. Thus the existence of a family of generalized convex
polytopes $\{P(t)|\: t \in (0,1]\}$ that satisfies the properties
\ref{ittt:1}.--\ref{ittt:3}. from Theorem \ref{thm:Family} is
established. It remains to show that $P(t)$ converges to a convex
polytope $P$ with a marked point $a$ as described in the property
\ref{ittt:5}.

\begin{lem} \label{lem:limr}
Let $(t_n)$ be a sequence in $(0,1]$ with $\lim_{n \to \infty} t_n
= 0$ such that there exist the limits
\begin{equation} \label{eqn:limr}
r_i(0) = \lim_{n \to \infty} r_i(t_n) \quad \mbox{for all }i.
\end{equation}
Then there exists a convex polytope $P \subset \R^3$ with boundary
isometric to $M$ and a point $a \in P$ such that $r_i(0) = \|p_i -
a\|$, where $p_i$ are the vertices of $P$.
\end{lem}
\begin{proof}. The lemma is obvious if $r(0) \in \P(M)$. Then $r(0)$ defines a generalized polytope $P(0)$ with curvatures
$$
\kappa_i(0) = \lim_{n \to \infty} \kappa_i(t_n) = 0.
$$
Therefore $P(0)$ is isometric to a convex polytope $P \subset
\R^3$ with a marked interior point $a$.

In the general case the following argument works. We may assume
that all triangulations $T(t_n)$ are equal to some triangulation
$T$. For every edge $ij \in \E(T)$ consider the dihedral angle
$\theta_{ij}(t)$ of the polytope $P(t)$ at $ij$. Angles
$\theta_{ij}(t)$ vary in $(0, \pi]$, so after choosing a
subsequence, if needed, we have the limits
$$
\lim_{n \to \infty} \theta_{ij}(t_n) = \theta_{ij} \in [0,\pi].
$$
Now pick up a face $F \in \F(T)$ and place it in $\R^3$. The radii
$r_i(0)$ for $i \in F$ determine the position of the point $a$, up
to the choice of an orientation. The dihedral angles at the edges
of $F$ determine the positions of the faces adjacent to $F$.
Besides, the condition $r_i(0) = \|p_i - a\|$ is satisfied for the
new vertices. The new faces determine the positions of their
neighbors, and so on. The faces fit nicely and bound a convex
polytope, if they fit around every vertex. Consider the sequence
of spherical sections $S^i_n$ at the vertex $i$ in the generalized
polytopes $P(t_n)$. This is a sequence of singular star polygons
with the angles $\theta_{ij}(t_n)$, curvature $\kappa_i(t_n)$ and
fixed edge lengths. It is not hard to prove that $S^i_n$ converges
to a convex spherical polygon with angles $\theta_{ij}$, possibly
degenerated to a doubly-covered spherical arc. This means that the
faces around the vertex $i$ can be put together. \end{proof}

Note that we just proved the existence part of Alexandrov's
theorem. However, properties \ref{ittt:4}. and \ref{ittt:5}. are
still to be proved.

Due to Lemma \ref{lem:Bound}, the radii $r_i(t)$ are bounded on
$(0,1]$. Therefore there is a sequence $t_n \to 0$ such that the
limits (\ref{eqn:limr}) exist. By Lemma \ref{lem:limr}, this gives
a convex polytope $P$ and a point $a \in P$. By the uniqueness
part of Alexandrov's theorem, the polytope $P$ does not depend on
the choice of a sequence $(t_n)$. It remains to prove that the
point $a$ is also independent of this choice. We will show that
the point $a$ satisfies the condition (\ref{eqn:apex}) and that
the condition (\ref{eqn:apex}) defines the point uniquely. This
will imply the properties \ref{ittt:4}. and \ref{ittt:5}. and thus
will complete the proof of Theorem \ref{thm:Family}.

Let $S_n$ be the spherical section of the generalized polytope
$P(t_n)$, see Definition \ref{dfn:SphSec}.
\begin{lem} \label{lem:Convergence}
Suppose that $S_n$ converges to a spherical section $S$ with the
vertices $v_i$. Then the following equality holds:
\begin{equation} \label{eqn:kv}
\sum_i \kappa_i(1) v_i = 0.
\end{equation}
\end{lem}
Some explanations are necessary. The spherical section $S_n$
inherits from the generalized polytope $P(t_n)$ the triangulation
$T(t_n)$. We say that the sequence of spherical sections $S_n$
converges, if $T(t_n) = T$ for almost all $n$, and for every
triangle $\Delta \in \F(T)$ the angles and the side lengths of its
spherical image $S\Delta_n$ converge. If $(t_n)$ is as in Lemma
\ref{lem:limr}, and the point $a$ is in the interior of $P$, then
the sequence $S_n$ converges to the spherical section $S$ of the
polytope $P$ viewed as a generalized polytope with the apex $a$.
If $a$ is an interior point of a face of $P$, then $S_n$ converges
to a spherical section $S$ that contains one degenerate triangle
with all angles equal to $\pi$. If $a$ is an interior point of an
edge, then $S_n$ does not necessarily converge, but can be
replaced by a converging subsequence. The limit section $S$
contains two degenerate triangles. In each of these cases the
vertices of the section $S$ are the unit vectors directed from $a$
to $p_i$:
$$
v_i = \frac{p_i - a}{\|p_i - a\|}.
$$
Finally, suppose $a = p_i$. Then $S_n$ contains a converging
subsequence. Its limit $S$ may contain degenerate triangles in the
star of $i$. However, $S$ has all curvatures zero and defines a
unit vector $v_i$, ``the direction from $a$ to $p_i$''.

\begin{proof} of Lemma \ref{lem:Convergence}. In the spherical section $S_n$ mark one of the singularities by $1$ and join it by geodesic arcs to all of the other singularities. The geodesics will not pass through singularities since they have positive curvature. Mark the other singularities by $2,\ldots,m$ in the clockwise order around $1$. Cut $S_n$ along the geodesic arcs and develop the result onto the unit sphere $\Sph^2$. Denote by $v_i(n)$ the images of the singularities $i, 2 \le i \le m$, and by $v_1(n)$ the image of $1$ that lies between $m$ and $2$. See Fig. \ref{fig:fig20}, where we omitted the index $n$. Clearly, by suitable rotations of $\Sph^2$ we can achieve $\lim_{n \to \infty} v_i(n) = v_i$.

\begin{figure}[ht]
\centerline{\input{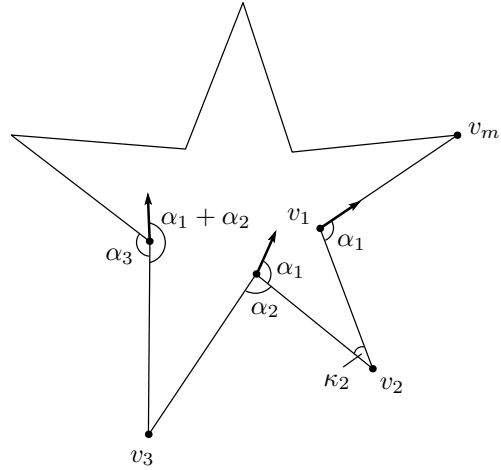}} \caption{Spherical
section $S_n$ developed onto the unit sphere. The development is
the exterior of the star.} \label{fig:fig20}
\end{figure}

Let $R_i \in SO(3)$ be the rotation around the vector $v_i(n)$ by
the angle $\kappa_i(t_n)$. We claim that
\begin{equation} \label{eqn:rot}
R_{m} \circ \cdots \circ R_2 \circ R_1 = \mathrm{id}.
\end{equation}
Indeed, from Fig. \ref{fig:fig20} it is clear that $R_{m} \circ
\cdots \circ R_2(v_1(n)) = v_1(n)$. On the same figure the
trajectory of a tangent vector at the point $v_1(n)$ is shown. The
vector gets rotated by the angle $\alpha_1(n) + \cdots +
\alpha_{m-1}(n) = \omega_1(n) = 2\pi - \kappa_1(n)$. The equality
(\ref{eqn:rot}) follows.

Since $\kappa_i(t_n) = t_n \cdot \kappa_i(1)$, we have for any
vector $x \in \R^3$
$$
R_i(x) = x + t_n \cdot \kappa_i(1) [v_i(n), x] + o(t_n),
$$
where $[\cdot,\cdot]$ denotes the cross product. Then
(\ref{eqn:rot}) implies
\begin{eqnarray*}
x = R_{m} \circ \cdots \circ R_2 \circ R_1(x) & = & x + \sum_{i=1}^{m} t_n \cdot \kappa_i(1) [v_i(n),x] + o(t_n)\\
& = & x + t_n \cdot \left[ \sum_{i=1}^{m} \kappa_i(1) v_i(n),x
\right] + o(t_n).
\end{eqnarray*}
It follows
$$
\sum_{i=1}^{m} \kappa_i(1) v_i(n) = o(1),
$$
and since $v_i(n) \to v_i$, we have (\ref{eqn:kv}). \end{proof}

{\em Remark.} When $a$ is an interior point of the polytope $P$,
Lemma \ref{lem:Convergence} can be proved using Proposition
\ref{prp:KerDescr} and the fact that the matrix
$(\frac{\partial\kappa_i}{\partial r_j})$ is symmetric.

\begin{lem} \label{lem:aPoint}
Let $P \subset \R^3$ be a convex polytope with the vertices $p_1,
\ldots, p_m$. Let $\delta_i$ be the angular defect at the vertex
$p_i$. Then there exists a unique point $a = a(P)$ in the interior
of $P$ (in the relative interior, if $P$ is $2$-dimensional) such
that
\begin{equation} \label{eqn:aPoint}
\sum_{i=1}^m \delta_i \frac{p_i - a}{\|p_i - a\|} = 0.
\end{equation}
\end{lem}
\begin{proof}. Consider the function $f:P \to \R$,
$$
f(x) = \sum_{i=1}^m \delta_i \cdot \|p_i - x\|.
$$
It is differentiable everywhere except the vertices of $P$, and
$$
\mathrm{grad}f(x) = -\sum_{i=1}^m \delta_i \frac{p_i - x}{\|p_i -
x\|}.
$$
Thus we are looking for the critical points of the function $f$.
Let $a$ be the point of the global minimum of $f$:
$$
f(a) = \min_{x \in P} f(x).
$$
Let us show that $a$ is a (relative) interior point of $P$.
Indeed, if $a$ is a boundary point of $P$ and not a vertex, then
the vector $-\mathrm{grad} f$ at $a$ is directed inside $P$ (which
should by understood appropriately if $P$ is $2$-dimensional).
Thus $a$ cannot be even a local minimum. Assume now $a = p_i$.
Then it is easy to see that
$$
\delta_i \ge \left| \sum_{j \ne i} \delta_j \frac{p_j - p_i}{\|p_j
- p_i\|} \right|.
$$
But this contradicts Lemma \ref{lem:Iso}. The existence of a point
$a$ that satisfies (\ref{eqn:aPoint}) is proved. The uniqueness
follows from the concavity of the function $f$. \end{proof}

If the numbers $\kappa_i(1)$ are sufficiently close to $\delta_i$,
then Lemma \ref{lem:Iso} and hence Lemma \ref{lem:aPoint} hold
also with $\delta_i$ replaced by $\kappa_i(1)$. Due to Lemma
\ref{lem:Convergence} this implies that for any sequence $(t_n)$
such that $P(t_n)$ converges to a convex polytope $P$, the apex of
$P(t_n)$ converges to the point $a \in P$ characterized by
(\ref{eqn:apex}). Theorem \ref{thm:Family} is proved.

\section{Lemmas from spherical geometry} \label{sec:Lemmas}
Lemmas of this section deal with spherical cone structures on the
disk or on the sphere. We consider only spherical cone structures
with positive curvatures at singularities.
\begin{lem} \label{lem:SphBasic}
The angles $\alpha, \beta, \gamma$ of a spherical triangle as
functions of side lengths $a,b,c$ have the following partial
derivatives:
\begin{eqnarray}
\frac{\partial \alpha}{\partial a} & = & \;\;\: \frac{1}{\sin b \sin \gamma}, \label{eqn:aa}\\
\frac{\partial \alpha}{\partial b} & = & -\frac{\cot \gamma}{\sin
b}. \label{eqn:ab}
\end{eqnarray}
The angles $\alpha, \beta, \gamma$ of a Euclidean triangle as
functions of side lengths $a,b,c$ have the following partial
derivatives:
\begin{eqnarray*}
\frac{\partial \alpha}{\partial a} & = & \;\;\: \frac{1}{b \sin \gamma}, \\
\frac{\partial \alpha}{\partial b} & = & -\frac{\cot \gamma}{b}.
\end{eqnarray*}
\end{lem}
\begin{proof}.
Direct calculation using spherical and Euclidean cosine and sine
theorems.
\end{proof}

As a technical tool we will use \emph{merging of singularities}.
Consider two cone points $A$ and $B$ with curvatures $\alpha$ and
$\beta$ on a surface with spherical cone structure. Suppose that
there is a geodesic arc joining $A$ and $B$. Such an arc always
exists, if all cone points have positive curvature, and the
boundary of the surface is convex. Cut the surface along this arc
and paste in a singular digon that is formed by two copies of a
triangle $ABC$ with angles $\frac{\alpha}{2}$ and
$\frac{\beta}{2}$ at vertices $A$ and $B$, identified along the
sides $AC$ and $BC$. The points $A$ and $B$ become regular,
instead there appears a new cone point $C$. Note that the
curvature of $C$ is less than $\alpha + \beta$.

\begin{dfn} \label{dfn:SingSphPol}
A singular spherical polygon is a 2-disk equipped with a spherical
polyhedral metric with piecewise geodesic boundary. It is called
convex, if the angles at the boundary vertices don't exceed $\pi$.
A convex star polygon is a convex singular spherical polygon with
a unique singularity.
\end{dfn}
The area of a singular spherical polygon can easily be computed:
\begin{equation} \label{eqn:SphArea}
\area = 2\pi - \sum_i \kappa_i - \sum_j (\pi - \theta_j),
\end{equation}
where $\kappa_i$'s are the curvatures of singularities, and
$\theta_j$ are the angles of the polygon, $\pi-\theta$ thus being
the outer angles.

\begin{figure}[ht]
\centerline{\input{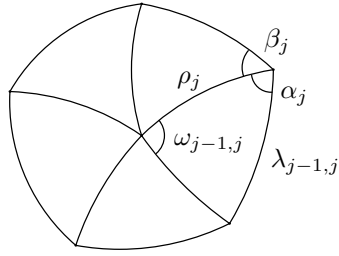}} \caption{Angles and
lengths in a convex star polygon.} \label{fig:fig12}
\end{figure}

A convex star polygon can be cut into a collection of spherical
triangles by a set of geodesic arcs joining the singularity to the
boundary. Or, other way round, a convex singular spherical polygon
can be specified as a result of gluing spherical triangles around
a common vertex. On Fig. \ref{fig:fig12} we introduce notation for
arc lengths and angles of a star polygon defined by a gluing
procedure. We allow $\alpha_j + \beta_j = \pi$ for some $j$, also
we allow a polygon be glued from a single triangle, two sides of
the triangle identified. Denote
$$
\kappa = 2\pi - \sum_j \omega_{j-1,j}
$$
the curvature of the singularity, and
$$
\delta = 2\pi - \sum_j \lambda_{j-1,j} = 2\pi - \per,
$$
where $\per$ stands for perimeter.

\begin{lem} \label{lem:Sph1}
If in a convex star polygon $\rho_j < \frac{\pi}{2}$ holds for all
$j$, then $\kappa < \delta$.
\end{lem}
\begin{proof}.
Consider the angles as functions of arc lengts. By formula
(\ref{eqn:ab}),
$$
\frac{\partial\kappa}{\partial\rho_j} = \frac{\cot\alpha_j +
\cot\beta_j}{\sin\rho_j} \ge 0.
$$
It is not hard to show that by increasing $\rho_j$'s and leaving
boundary edges constant we can deform any polygon with $\rho_j <
\frac{\pi}{2}$ into a polygon with $\rho_j = \frac{\pi}{2}$ for
all $j$. If $\rho_j = \frac{\pi}{2}$ for all $j$, then $\kappa =
\delta$. Hence the statement of the Lemma.
\end{proof}

\begin{lem} \label{lem:Sph2}
If in a convex star polygon $P$ holds $0 < \kappa < \delta$, then
its boundary cannot contain a geodesic arc of length $\pi$ or
more.
\end{lem}
\begin{proof}.
Let us assume the converse. Choose on the boundary of $P$ two
points $A$ and $B$ that divide the boundary in two curves $L$ and
$L'$ such that $L$ is geodesic and has length $|L| = \pi$. Cut $P$
in two spherical polygons $C$ and $C'$ by geodesic arcs $OA$ and
$OB$, where $O$ is the singularity of $P$ and $C$ contains $L$ in
its boundary. Clearly, the polygon $C$ is a spherical lune. Thus
it has angle $\pi$ at the vertex $O$. Denote by $\gamma$ the angle
of the polygon $C'$ at the vertex $O$. We have $\gamma = \pi -
\kappa < \pi$. Thus $C'$ is convex and contains a geodesic arc
$AB$. Consider the triangle $ABO$ contained in $C'$, and denote
its side lengths $a,b$, and $c$.

\begin{figure}[ht]
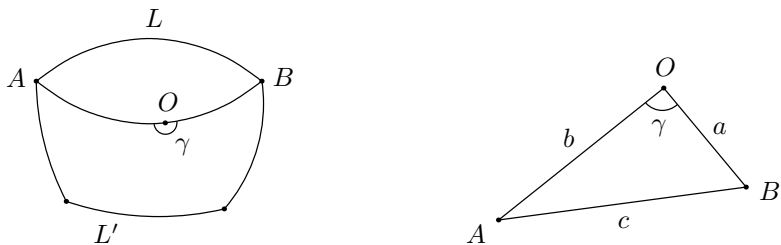

\centerline{\input{Figures/fig13a.pstex_t}\hspace{2.5cm}\input{Figures/fig13b.pstex_t}}
\caption{To the proof of Lemma \ref{lem:Sph2}.} \label{fig:fig13}
\end{figure}

We have
\begin{eqnarray*}
\kappa & = & \pi - \gamma, \\
\delta & = & \pi - |L'| \le \pi - c.
\end{eqnarray*}
Also $a + b = \pi$. By the spherical cosine theorem,
\begin{eqnarray*}
\cos c & = & \cos a \cos b + \sin a \sin b \cos \gamma \\
       & = & -\cos^2 a + \sin^2 a \cos \gamma.
\end{eqnarray*}
It follows
$$
1+ \cos c = \sin^2 a (1+\cos \gamma) < 1 + \cos \gamma.
$$
Thus $c > \gamma$ and $\kappa > \delta$ which is a contradiction.
\end{proof}

%Consider the sector $S$ based on a boundary arc of length $\pi$. The rest of the polygon is a non-singular spherical polygon $P$. Denote the angles as shown on Fig. \ref{fig:Sph2}. First note that $\gamma$ is less than $\pi$, otherwise the perimeter of the initial polygon is larger than $\pi + a + b > 2\pi$ that contradicts the assumption $\delta > 0$. Thus $P$ is a convex spherical polygon, and we can join the endpoints of the long boundary arc by an arc inside $P$. This gives us a triangle $\Delta$ pictured in the right part of Fig. \ref{fig:Sph2}.

%Let us show that $\gamma' \le \pi$. The contrary assumed, we have $\alpha' + \beta' > \pi$. But $\alpha + \beta \ge \pi$ (since $a+b \ge \pi$, the triangle $\Delta$ constitutes at least a half of the lune with angle $\gamma$, hence $\alpha + \beta + \gamma - \pi \ge \gamma$). Thus either $\alpha + \alpha' > \pi$ or $\beta + \beta' > \pi$ which contradicts to convexity of the initial polygon.

%Compare the triangle $\Delta$ with a triangle $\Delta'$ which has the angle $\gamma'$ between sides $a$ and $b$. Triangle $\Delta'$ complements the sector $S$ to a half-sphere. Therefore its side opposite to $\gamma'$ has length $c' = \pi$. Deform $\Delta$ to $\Delta'$ keeping $a$ and $b$ fixed. Because of $\frac{\partial\gamma}{\partial c} = \frac{\csc\beta}{\sin a} > 1$ we must have $\gamma' - \gamma > c' - c$. But $\gamma' - \gamma = \kappa$ and $\per \ge \pi + c$ lead to a contradiction with $\kappa < \delta$.

\begin{lem} \label{lem:Sph3}
For any $0 < c_1 < c_2 < \pi$ and $\gamma > 0$ there exists an
$\epsilon > 0$ such that the following holds. Any convex singular
spherical polygon with all side lengths between $c_1$ and $c_2$
and all outer angles greater than $\gamma$ has perimeter less than
$2\pi - \epsilon$.
\end{lem}
In other words, a polygon with outer angles not too small and
sides neither too large nor too small cannot have perimeter too
close to $2\pi$. Note that the lemma is false if the polygon is
allowed to have singularities of negative curvature.
\begin{proof}.
First consider the case when the polygon is non-singular. Let $a$
be the length of one of its sides. Then the polygon is contained
in the spherical triangle with one side of length $a$ and both
adjacent angles $\pi - \gamma$; hence its perimeter is less than
that of this triangle. For $a$ between $c_1$ and $c_2$ the
perimeters of corresponding triangles are uniformly bounded from
above by $2\pi - \epsilon$ for some $\epsilon > 0$. Thus in the
non-singular case assumptions of the Lemma can be weakened:
inequalities for one side and two adjacent angles already suffice.

\begin{figure}[ht]
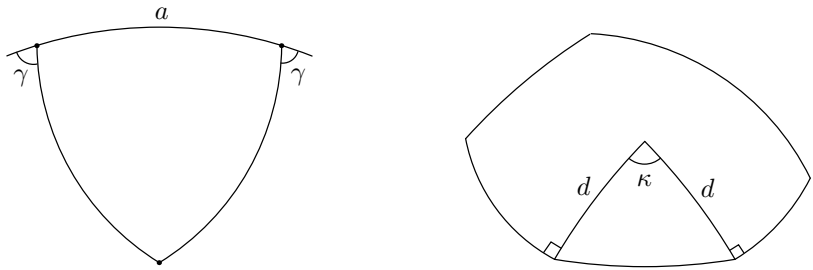

\centerline{\input{Figures/fig14a.pstex_t}\hspace{2cm}\input{Figures/fig14b.pstex_t}}
\caption{To the proof of Lemma \ref{lem:Sph3}.} \label{fig:fig14}
\end{figure}

If the polygon is singular, merge all of its singularities
consequtively. Let $\kappa$ be the curvature of the resulting star
polygon $P$. The merging changes neither the sides nor the angles
of the polygon, so it suffices to prove the Lemma for the polygon
$P$. If $P$ has only one or two sides, then its perimeter is less
than $2c_2$. Otherwise draw the shortest arc joining the
singularity to the boundary of the polygon $P$; let $d$ be its
length. We claim that $d < \frac{\pi}{2}$. Indeed, embed $P$ into
the sphere with two singularities of curvature $\kappa$. This
sphere is obtained from the usual sphere by removing a lune of
angle $\kappa$ and identifying the boundary semicircles. If $d \ge
\frac{\pi}{2}$, then the $\kappa$-sphere can be covered by two
copies of $P$. Hence the area of $P$ has to be at least $2\pi -
\kappa$. But it is smaller than this by formula
(\ref{eqn:SphArea}). Thus $d < \frac{\pi}{2}$. Cut along the
shortest arc and glue to the sides of the cut an isosceles
triangle with the angle $\kappa$ at the vertex. Due to $d <
\frac{\pi}{2}$ the angles at its base are less than
$\frac{\pi}{2}$. Thus the resulting non-singular polygon is
convex. Besides, at least one of its sides together with the
adjacent angles satisfy assumptions of the Lemma. By the previous
paragraph, the perimeter of the resulting polygon is less than
$2\pi - \epsilon$. Then so is the perimeter of the initial
singular polygon.
\end{proof}

\begin{lem} \label{lem:Sph4}
Let a spherical polyhedral metric on the sphere be given, with all
of the singularities of positive curvature. Let one of the
singularities be labeled by $0$. Then
$$
\kappa_0 \le \sum_{i \ne 0} \kappa_i.
$$
Suppose additionally that the singularity $0$ is the center of a
star polygon $C$ with vertices at some of the other singularities
and all angles at least $\pi$. Let $\per(C)$ denote the perimeter
of $C$. Then the following inequality holds:
\begin{equation} \label{eqn:Last}
\kappa_0 \ge \left( 1 - \frac{\per(C)}{2\pi} \right) \, \sum_{i
\ne 0} \kappa_i.
\end{equation}
\end{lem}
\begin{proof}.
To prove the first inequality, merge consequtively all the
singularities except of $0$. The result is the sphere with two
singularities of curvature $\kappa_0$. Since merging decreases the
total curvature, the result follows.

The proof of the second part of the Lemma is much more involved.
Denote by $D$ the convex singular polygon that is the complement
to the star polygon $C$. We make $D$ non-singular, increasing at
the same time the right hand side of (\ref{eqn:Last}). This is
done by a sequence of operations inverse to merging of
singularities. Choose a vertex of $D$, denote it by $1$ and join
it to a singularity $i$ in the interior of $D$ by the shortest arc
$a_i$. Extend the arc $a_i$ beyond $i$ by a geodesic $b_i$ such
that the two angles between $a_i$ and $b_i$ are equal. The arc
$b_i$ ends either at some other singularity $j$ or on the boundary
of $D$. If the latter is the case and if $D$ contains two
spherical triangles with the sides $a_i$ and $b_i$, then say that
we have good luck. In this case we cut both triangles out and glue
their remaining sides together, see Figure \ref{fig:fig15}. This
splitting of the singularity $i$ increases the right hand side of
(\ref{eqn:Last}), because it does not change the perimeter and
increases the sum $\sum_{i \ne 0} \kappa_i$. On the other hand it
decreases the number of singularities inside $D$, so that we can
proceed inductively. In the case of bad luck draw the arcs $a_j$
and $b_j$ for all singularities $j$ from the vertex $1$. Then the
arc $b_i$ crosses some $a_j$ or $b_j$. The former is not possible
because then $a_j$ is not the shortest. If $b_i$ ends at the
singularity $j$ without crossing any $b_k$ before this, then we
cut out a pair of equal triangles with vertices $i,j$ and $1$,
thus decreasing the number of singularities. Finally, if $b_i$
crosses the interior $b_j$ at the point $B_{ij}$, we can assume
without loss of generality that neither of the arcs $iB_{ij}$ and
$jB_{ij}$ is crossed by some $b_k$. Then we can cut out two pairs
of equal triangles with the common vertex $B_{ij}$ and again
reduce the number of singularities.

\begin{figure}[ht]
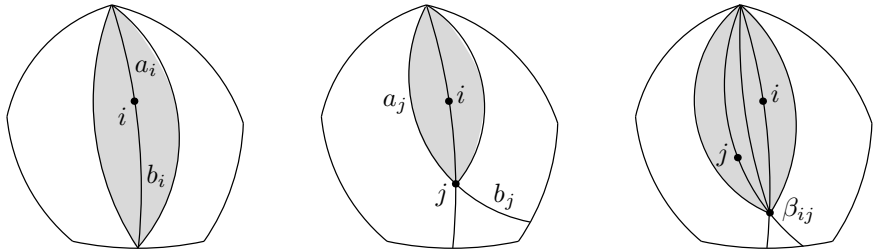

\centerline{\input{Figures/fig15a.pstex_t}\hspace{1cm}\input{Figures/fig15b.pstex_t}\hspace{1cm}\input{Figures/fig15c.pstex_t}}
\caption{Splitting the singularities inside the singular polygon
$D$.} \label{fig:fig15}
\end{figure}

Now the Lemma is reduced to the case of non-singular $D$. Embed
the star polygon $C$ into the sphere with two singularities of
curvature $\kappa_0$ and denote by $D'$ its complement. Denote by
$\theta_i$ and $\theta'_i$ the respective angles of $D$ and $D'$.
We have $\theta'_i - \theta_i = \kappa_i > 0$, whereas the
respective sides of $D$ and $D'$ are equal, see Figure
\ref{fig:fig16}. Therefore the inequality (\ref{eqn:Last}) can be
viewed as a comparison inequality between a usual spherical
polygon $D$ and a star polygon $D'$, both $D$ and $D'$ convex.

\begin{figure}[ht]
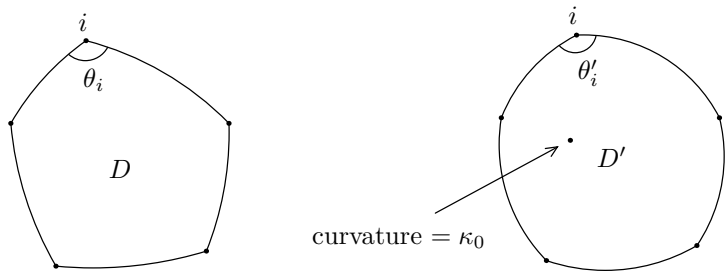

\centerline{\input{Figures/fig16a.pstex_t}\hspace{1cm}\input{Figures/fig16b.pstex_t}}
\caption{A convex polygon $D$ and a convex singular polygon $D'$;
$\theta'_i - \theta_i = \kappa_i$.} \label{fig:fig16}
\end{figure}

Note that the star polygon $D'$ is completely determined by the
polygon $D$ and angle differences $\kappa_i$. To see this, cut
$D'$ along a geodesic joining the central singularity with any
vertex and develop the result on the sphere. On the other hand, a
procedure similar to the merging of singularities allows us to
obtain $D'$ from $D$ as follows. Cut $D$ along a diagonal $ij$ and
glue in a digon with the angles $\kappa_i$ and $\kappa_j$. The
result is a star polygon with the angles $\theta'_i$ and
$\theta'_j$ at the vertices $i$ and $j$, respectively. Then glue
in a digon between a boundary vertex $k$ and the central
singularity so that the singularity disappears and the angle at
$k$ becomes $\theta'_k$. This yields a new convex star polygon,
and so on.

If $\theta_i = \theta'_i$ for some $i$, then the above
transformation of $D$ into $D'$ does not involve the vertex $i$ at
all. If $D$ and $D'$ have more than two sides, this implies that
the polygon $D'$ contains a non-singular triangle $\Delta'_i$
spanned by the sides incident to the vertex $i$. The triangle
$\Delta'_i$ is congruent to the corresponding triangle $\Delta_i$
in $D$. If we cut both triangles off, then we get a new couple of
polygons with the same $\kappa_0$ and $\sum_{i \ne 0} \kappa_i$
but with a smaller perimeter. Thus the inequality (\ref{eqn:Last})
for the new couple implies (\ref{eqn:Last}) for the old.

The rest of the proof goes by induction on the number of sides: we
carefully deform $D$ until $\theta_i$ becomes equal to $\theta'_i$
for some $i$, and cut off the triangles $\Delta$ and $\Delta'$ as
in the previous paragraph. The induction base is two sides: $D$ is
a doubly covered segment, $D'$ is a symmetric star digon. Let us
deal with the induction step first. Suppose that $D$ has three or
more sides, and $\theta'_i > \theta_i$ for all $i$. Choose a
vertex $i$ and denote by $a$ and $b$ the lengths of the adjacent
sides. Denote by $c$ the length of the diagonal in $D$ (or side,
if $D$ is the triangle) that forms a triangle together with these
two sides. Without loss of generality, $a \le b$. Consider the
following deformation of $D$. All of the vertices except of $i$
stay fixed, and $i$ moves so that the sum $a+b$ remains constant
and $a$ decreases. Define the corresponding deformation of $D'$ as
follows. Develop $D'$ onto the sphere by cutting it along a radius
that ends at a vertex other than $i$. The deformation leaves all
of the vertices of the development fixed, except of $i$ that moves
so that $a+b$ remains constant and $a$ decreases, as in $D$.
Denote by $c'$ the length of the diagonal that corresponds to the
diagonal $c$ in $D$. See Figure \ref{fig:fig17}.

\begin{figure}[ht]
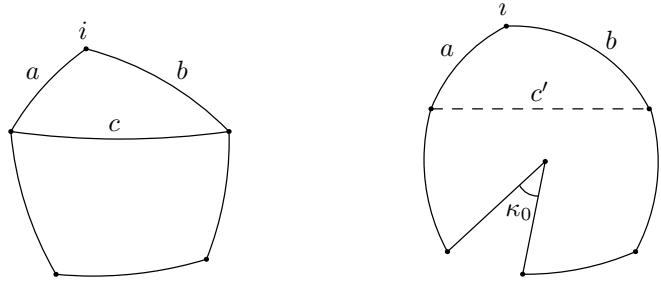

\centerline{\input{Figures/fig17a.pstex_t}\hspace{2.5cm}\input{Figures/fig17b.pstex_t}}
\caption{Polygon $D$ and singular polygon $D'$ cut along a
radius.} \label{fig:fig17}
\end{figure}

Clearly, the deformation of $D'$ does not change $\kappa_0$. By
construction, the perimeter of $D$ and $D'$ is also constant. As
we will show, the deformation increases the sum $\sum \kappa_i$.
This implies that during the deformation the left hand side of
(\ref{eqn:Last}) remains constant, whereas the right hand side
increases. The deformation stops in one of the following cases.

\emph{We get $\theta_j = \theta'_j$ for some $j$.} Then we cut off
the triangles $\Delta$ and $\Delta'$ at the vertex $j$ and use the
induction assumption.

\emph{We get $\theta'_j = \pi$ for some $j$.} In this case we
start the deformation that affects the vertex $j$. It does not
change the angles of $D'$, so it has to end for a different
reason.

\emph{The polygon $D$ is a triangle and it degenerates so that $b
= a+c$.} In this case $D'$ degenerates to a symmetric digon. This
is the induction base.

\emph{The side $b$ hits the image of the singularity in the
development of $D'$.} This is impossible, because if the
singularity is close to the boundary of $D'$ but far from its
vertices, then $\kappa_0$ is small. This follows from the
transformation of $D$ into $D'$ by gluing in digons as described
before.

In order to show that $\sum \kappa_i$ increases during the
deformation of $D$ and $D'$, consider the triangle with side
lengths $a,b,c$ and the triangle with side lengths $a,b,c'$. It is
not hard to see that
$$
d \sum \kappa_i = d \area' - d \area,
$$
where $\area$ and $\area'$ are the areas of the first and of the
second triangle, respectively. If $\alpha, \beta, \gamma$ are the
angles of the first triangle, then we have
$$
d \area = \frac{\partial \area}{\partial b} - \frac{\partial
\area}{\partial a} = \left( \frac{\partial}{\partial b} -
\frac{\partial}{\partial a} \right)(\alpha + \beta + \gamma).
$$
Since $c' > c$, for $d \area' > d \area$ it suffices to show
$\frac{\partial}{\partial c}d \area > 0$. With the help of Lemma
\ref{lem:SphBasic} one can compute
$$
\frac{\partial^2 \area}{\partial b \partial c} = \frac{\sin\alpha
- \sin(\beta + \gamma)}{\sin^2\alpha \sin\beta \sin\gamma \sin b
\sin c}.
$$
The denominator is invariant under permutations of $a,b,$ and $c$.
Thus we have
\begin{eqnarray*}
d \sum \kappa_i & = & \frac{\sin\alpha - \sin\beta + \sin(\alpha + \gamma)
- \sin(\beta + \gamma)}{\sin^2\alpha \sin\beta \sin\gamma \sin b \sin c} \\
& = & \frac{2\sin{\frac{\alpha - \beta}{2}} \cos{\frac{\alpha +
\beta + \gamma}{2}} \cos{\frac{\gamma}{2}}}{\sin^2\alpha \sin\beta
\sin\gamma \sin b \sin c}.
\end{eqnarray*}
This is easily seen to be positive.

Now let us prove the induction base. The digon $D'$ is glued from
two equal triangles. If we use the standard notation for the sides
and angles of a triangle, the inequality (\ref{eqn:Last}) becomes
\begin{equation} \label{eqn:theverylast}
\pi - \gamma > \left( 1 - \frac{c}{\pi} \right)(\alpha + \beta),
\end{equation}
where $\alpha, \beta \le \frac{\pi}{2}$ since $D'$ is convex. Let
us show that for $c$ and $\alpha + \beta$ constant, the left hand
side of (\ref{eqn:theverylast}) attains its minimum at $\alpha =
\beta$. Assume $\alpha < \beta$. Then $a < b$. It follows that the
simultaneous increase of $\alpha$ and decrease of $\beta$ by the
same amounts result in an increase of the area of the triangle.
Hence the left hand side of (\ref{eqn:theverylast}) decreases. For
$\alpha = \beta$, we have to show
\begin{equation} \label{eqn:final}
\frac{\pi - \gamma}{\pi - c} > \frac{2\alpha}{\pi},
\end{equation}
for an isosceles triangle based on $c$. We have $\cos
\frac{\gamma}{2} = \sin \alpha \cos \frac{c}{2}$. For $\alpha$
constant, the left hand side of (\ref{eqn:final}) can be shown to
be a monotonically increasing function of $c$. As $c \to 0$, it
tends to $\frac{2\alpha}{\pi}$, which proves the Lemma.
\end{proof}

\begin{lem} \label{lem:Iso}
Let $P \subset \R^3$ be a convex polytope with vertices $(p_i)$
and angle defects $(\delta_i)$. Then for any $i$
$$
\left| \sum_{j \ne i} \delta_j \frac{p_j - p_i}{\|p_j - p_i\|}
\right| > \delta_i.
$$
\end{lem}
\begin{proof}. Let $C \subset \Sph^2$ be the spherical section of the polytope $P$ at the vertex $i$. For any $j \ne i$ we have
$$
v_j = \frac{p_j - p_i}{\|p_j - p_i\|} \in C.
$$
The perimeter of $C$ equals $2\pi - \delta_i$. By Lemma
\ref{lem:c}, there exists a point $c \in C$ such that
$$
\dist(v_j,c) \le \frac{\pi}{2} - \frac{\delta_i}{4} \quad
\mbox{for all } j.
$$
Thus we have
$$
\left\langle \sum_{j \ne i} \delta_j v_j, c \right\rangle =
\sum_{j \ne i} \delta_j \langle v_j, c \rangle \ge \sum_{j \ne i}
\delta_j \cdot \sin\frac{\delta_i}{4} = (4\pi - \delta_i)
\sin\frac{\delta_i}{4} > \delta_i,
$$
where the last inequality follows easily from $\sin x >
\frac{2x}{\pi}$ for $0 < x < \frac{\pi}{2}$. \end{proof}

\begin{lem} \label{lem:c}
Let $C$ be a convex spherical polygon of perimeter $\per$. Then
$C$ is contained in a circle of radius $\frac{\per}{4}$. That is,
there exists a point $c \in \Sph^2$ such that
$$
\dist(x,c) \le \frac{\per}{4} \quad \mbox{for all } x \in C
$$
in the intrinsic metric of the sphere.
\end{lem}
\begin{proof}. Let $O$ be the circle of the smallest radius $r$ that contains $C$. Then there are three vertices $v_1, v_2$ and $v_3$ of $C$ such that $v_1, v_2, v_3$ lie on $O$, and the center of $O$ lies in the triangle $\Delta$ with the vertices $v_1, v_2, v_3$. Since $\per(\Delta) \le \per$, it suffices to prove
\begin{equation} \label{eqn:Triangle}
r \le \frac{\per(\Delta)}{4}.
\end{equation}
The triangle $\Delta$ is defined by the radius $r$ and the angles
$\alpha,\beta,\gamma \in (0,\pi)$ between the radii drawn from the
center of $O$ to the vertices. By computing the derivatives
$\frac{\partial \per(\Delta)}{\partial \alpha} - \frac{\partial
\per(\Delta)}{\partial \beta}$ one sees that the minimum of
$\per(\Delta)$ is achieved when $\Delta$ degenerates to an arc. In
this case we have an equality in (\ref{eqn:Triangle}). The lemma
is proved. \end{proof}

\bibliographystyle{abbrv}
\bibliography{Alexandrov}

\end{document}